\documentclass[reqno,12pt]{amsart}
\usepackage{amsmath}
\usepackage{amsfonts}
\usepackage{amssymb}
\usepackage{graphicx}
\usepackage{cite}
\usepackage[usenames]{color}
\usepackage{colortbl}
\usepackage[latin1]{inputenc}
\usepackage{hyperref}
\usepackage{esint}

 \newcommand{\R}{\mathbb R}
 \newcommand{\N}{\mathbb N}
 \newcommand{\dis}{\displaystyle}
 \renewcommand{\d}{\mathrm d}

\newcommand{\beq}{\begin{equation}}
\newcommand{\eeq}{\end{equation}}

\newcommand{\md}{\medskip}

\newtheorem{thm}{Theorem}
\newtheorem{lem}[thm]{Lemma}
\newtheorem{prop}[thm]{Proposition}
\newtheorem{cor}[thm]{Corollary}
\newtheorem{rem}{\it Remark}

\newtheorem{defin}{\it Definition}
\newenvironment{pf}{\it Proof. \enskip \rm}{\hfill $\qed$}

\begin{document}

\title[]{Finite propagation and saturation 
in reaction-diffusion-advection equations \\ governed by p-Laplacian operator}
\author{ Cristina Marcelli}

\address[Cristina Marcelli]{Dipartimento di Ingegneria Industriale e Scienze Matematiche,
Universit\`a Politecnica delle Marche, Via Brecce Bianche 12, 60131 Ancona, Italy}
\email{c.marcelli@staff.univpm.it}


\date{}

\maketitle

\begin{small}
\noindent {\bf Abstract.} 
The paper concerns front propagation for the following mono-stable reaction-diffusion-advection equation
\[f(u)u_x + g(u)u_\tau = [d(u)|u_x|^{p-2} u_x]_x+ \rho(u), \quad 
(x,\tau)\in \R\times [0,+\infty).\]
Besides existence and non-existence results for traveling wave solutions, the main focus is their classification: we provide criteria to establish if they attain one or both the equilibria at a finite time and in this case, if they are continuable as $C^1$-solutions or if they are sharp solutions.
\end{small}

\bigskip

\section{Introduction}

This paper concerns the existence and the asymptotic properties of the traveling wave solutions (t.w.s.) of the following general reaction-diffusion-advection equation governed by the p-Laplacian differential operator:
\beq \label{eq:RD}
 f(u)u_x + g(u)u_\tau = [d(u)|u_x|^{p-2} u_x]_x+ \rho(u), \quad 
(x,\tau)\in \R\times [0,+\infty) \eeq
where $p>1$, $f,g,\rho \in C[0,1]$, $d\in C^1(0,1)$.

 We deal with the so-called mono-stable equations, that is the reaction term $\rho$ is assumed to satisfy
\[
\rho(u)>0 \ \text{ in } (0,1), \quad \rho(0)=\rho(1)=0. \tag{$H_\rho$}
\]
Regarding the advection term f, we do not require any further condition other than continuity. Instead, as for the 
 the accumulation term $g$,  we assume that 
\[
\label{ip:g}
g(0)>0  \quad \text{  and } \quad  \int_0^u g(s)  \d s>0 \ \text{ for every }  u\in (0,1]. \tag{$H_g$}
\]
Finally, 
as for the diffusion term $d$, we suppose

\[
d(u)>0 \ \text{ in } (0,1) 
%\quad  \lim_{u\to 0} \rho(u)(d(u))^\frac{1}{p-1} =0  
\quad  \text{ and } \quad \rho d^{\frac{1}{p-1}} \in L^1(0,1). \ \tag{$H_d$}
\]
Notice that we do not assume $d$ is continuous at the equilibria, nor it is locally bounded.
\smallskip 
Throughout the paper we will always assume that the structural  hypotheses $(H_\rho)$, $(H_g)$ and $(H_d)$  listed above are satisfied, without needing to recall them each time.

\medskip
As it is well known, a t.w.s. for equation \eqref{eq:RD} is a solution of the type $u(x-c\tau)$ for some constant $c\in \R$, called  {\em wave speed}. So, the right t.w.s. connecting the equilibria 0 and 1 are solutions of the following boundary vale problem

\beq \label{odeprobold} \begin{cases} (d(u)|u'|^{p-2}u')' + (c g(u)-f(u)) u' +\rho(u)=0 \\ u(\alpha)=1 \ , \ \ u(\beta)=0, \end{cases} \eeq
in some interval $(\alpha,\beta)\subseteq \R$ (the left t.w.s. can be obtained by reflection and satisfy the reverse boundary conditions). 
In the monostable case, it is well known that there exists a threshold wave speed $c^*$ such that the dynamic supports t.w.s. if and only if $c\ge c^*$.
There is extensive literature on the existence of traveling wave solutions and on estimates of the minimal admissible wave speed
 $c^*$. We limit ourselves to quote some results for equations governed by the $p-$Laplacian operator, see  e.g. \cite{AV1,GS2,Aud1,AV2}  and \cite{DZ,DJKZ} for existence results in the case of discontinuous coefficients. Moreover, see \cite{cmp1} for equations involving also a non-constant accumulation term $g$.

The interest related to the existence or non-existence of t.w.s. and the admissible wave speeds $c$ naturally arises from the study of the behavior for large times of the solutions of \eqref{eq:RD},  which evolve as a t.w.s. (see e.g. \cite{MR}, \cite{DT1} or the recent paper \cite{AHR} and references therein contained). In this context, a relevant topic concerns the structure of the interval $(\alpha,\beta)$: if it is upper and/or lower bounded.  In fact, the question if a t.w.s. attains one or both the equilibria at finite times has a particular relevance in view of the applications. If $u(t_0)=1$ for some $t_0$
 and $u$ represents population density,  the model predicts that the population saturates the environment within a finite time, indicating  {\em finite saturation}. Instead, when $u(t_0)=0$ for some $t_0$ then this means that the area in which a portion of the population has spread always remains bounded, that is the dynamic is said to exhibits the phenomenon of {\em finite propagation} and describes the propagation, with a finite speed, of an interface that separates the region where $u>0$ from that where $u=0$. 

The occurrence of finite propagation is a well-known phenomenon in the case of linear differential operator ($p=2$) with density-dependent diffusion satisfying $d(0)=0$: in this case  the t.w.s. having the minimal wave speed $c^*$ reaches the steady state 0 at a finite time with a non-zero slope ({\em sharp} t.w.s.), while all the t.w.s. having wave speed $c>c^*$   approach the equilibrium with null slope ({\em front} t.w.s.): see, e.g., \cite{dPV}, \cite{GK},  \cite{MM1}, \cite{bcm2}, \cite{cmp1}. Finally, we mention the paper \cite{SRJ} which contains some numerical simulations showing the emergence of sharp-type t.w.s. at the equilibrium 0.
Nevertheless, notice that even when $c>c^*$ the t.w.s. can attain the equilibrium 0 at a finite time: indeed, in our general framework, in which the coefficients can be not Lischitz continuous and the equation can be degenerate or singular,  the  
uniqueness of the solution of equation in \eqref{odeprobold} satisfying $u(t_0)=u'(t_0)=0$ for some $t_0$ is not guaranteed, so equation \eqref{odeprobold} could admit non-constant regular solutions satisfying $u(t)=u'(t)=0$ for $t$ large enough. So, the finite propagation can occur even for $c>c^*$.

Summarizing, as for the behavior near the equilibium 0 we can have three types of solutions: {\em front-type} (i.e. regular) with $u(t)>0$ for every $t\in \R$, front-type with $u(t)=u'(t)=0$ for every $t\ge t_0$, and {\em sharp-type} that is with $u(t_0)=0$ and $u'(t_0^-)<0$ for some $t_0\in \R$ (see the figures below). 

\medskip
\includegraphics[scale=0.6]{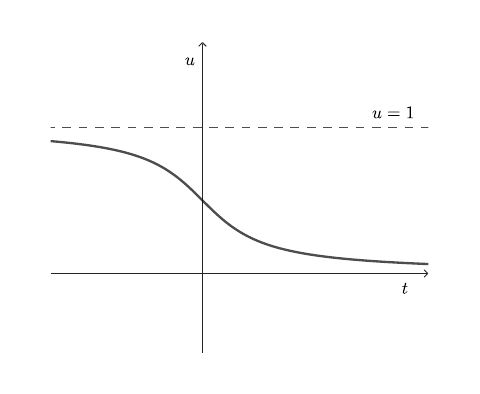}
\includegraphics[scale=0.6]{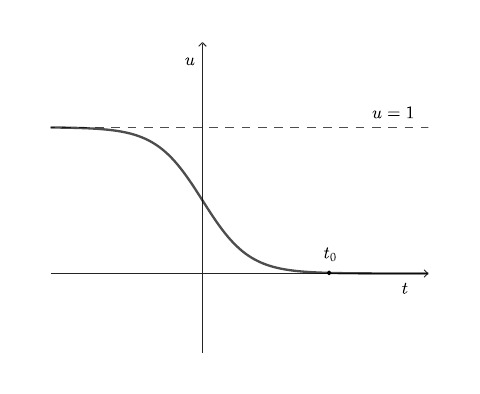}
\includegraphics[scale=0.6]{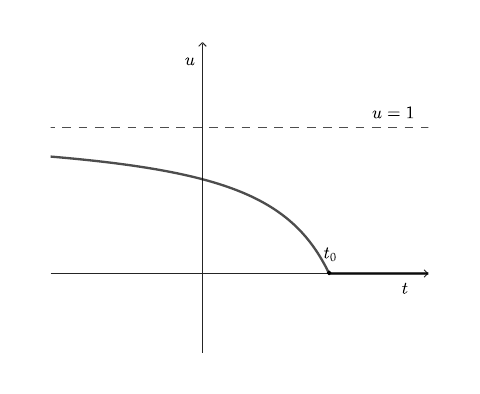}

While the classification between front and sharp solutions (that is between $C^1$ and non regular t.w.s) has been widely studied in the case $p=2$, the same does not apply to 
the distinction between the first two types of solutions, 
and therefore the possible appearance of finite propagation even with $C^1$-solutions is not yet fully understood.
 Drabek et al. carried out a first investigation for equations without convection  governed by the p-laplacian: in \cite{DZ} in the case of lack of convective effects and in \cite{DJKZ} for convective-reaction-diffusion equations. They showed that the finite propagation does not occur when $c$ is large enough, but the study for each $c>c^*$ remained an open problem. Moreover, also the behavior of the t.w.s. having speed $c^*$ in the case of equations involving the $p$-Laplacian was not known yet.

\medskip
A similar discussion can be made also for the steady state $u=1$: we can have front-type solution such that $u(t)<1$ everywhere, or front-type solution attaining the value 1 at a finite time (with null slope) or sharp-type solutions. Hovewer, the behaviour of the solutions when approaching 1 (stable equilibrium) is non analogous to that occurring for 0 (unstable equilibrium) and the results are completely different.
In the case of linear differential operators, we can quote \cite{MM1}, \cite{bcm2} and \cite{cmp1} for criteria concerning the regularity of the solutions (distinguishing front-type from sharp-type), whereas in  \cite{DT1} the investigation possible attaining of the equilibrium 1 at finite times was carried out, for equations without advection terms. Finally, in \cite{DZ}, \cite{DJKZ} some criteria for the occurrence of finite saturation are presented in the case of the $p-$Laplacian differential operator, but also in this case, the results were not complete in the case of equation with convective effects (see Remark \ref{rem:DR}).

\medskip
In this paper we consider the general equation \eqref{eq:RD}, governed by the $p$-Laplacian operator and including both a convective term $f$ and an accumulation one $g$. 
Although recent studies address discontinuous coefficients (see
\cite{DZ} and \cite{DJKZ}), we assume continuity of all coefficients in (0,1), as our main focus is the local asymptotic analysis at equilibria 0 and 1, which is unaffected by jump discontinuities.

After providing an existence and non-existence result (see Theorem \ref{t:main2}) together with an estimate for $c^*$ (see also Theorem \ref{t:stima}), we develop the local analysis at both the equilibria, providing conditions, most of which necessary and sufficient, to classify the t.w.s. distinguishing between the three types of solutions. 

More in detail, 
as regards the finite propagation, we show that 
this occurs for $c>c^*$ if and only if the integral $\dis\int_0^\frac12 \frac{1}{\rho(u)} \d u$ converges (see Theorem \ref{t:T2}). Moreover, we show that the appearance of finite propagation in the case $c=c^*$ depends, besides the convergence of the previous integral, also on the value $p$ (see Theorem \ref{t:T2*}). In the particular case in which $d$ and $\rho$ are asymptotic to a power $u^\gamma$ as $u\to 0$, we  summarize the results in Corollary \ref{cor:0}. Finally, since the asymptotic behavior also depends on the value $c^*g(0)-f(0)$ (if it is positive or null), we also provide a result concerning this topic (see Theorem \ref{t:stima}).

\smallskip
Section \ref{sec:1} is devoted to the investigation of the possible appearance of the phenomenon of finite saturation. We provide some criteria, most of which necessary and sufficient, to classify the t.w.s. distinguishing between front and sharp (i.e. between $C^1$ and non regular), and among the front-type solutions, between those satisfying $u(t)<1$ everywhere from those attaining the value 1 at some point. More in particular, we show that also in this case, the asymptotic behavior depends on the sign of the value $cg(1)-f(1)$: we provide results both for the case $c>f(1)$ and for the case $c<f(1)$, which was not treated in \cite{DJKZ} (see Corollary \ref{cor:1}).

\bigskip

\section{T.w.s. and their properties}

First of all, we clarify what we mean by a t.w.s. of \eqref{eq:RD}. 

\begin{defin}\label{def:sol}
A traveling wave solution {\rm (}t.w.s.{\rm)} to  \eqref{eq:RD} is a function \(u\in C^1(\alpha,\beta)\), with \((\alpha,\beta)\subseteq \R\), such that \(0<u(t)<1\) for every \(u\in (\alpha,\beta)\), \((d\!\circ \!u)|u'|^{p-2}u'\in C^1(\alpha,\beta)\), such that
\beq \label{eq:ode_lapl} (d(u)|u'|^{p-2}u')' + (c g(u)-f(u)) u' +\rho(u)=0 \quad \text{ for every } t \in (\alpha,\beta);
\eeq
\beq \label{c:as} u(\alpha^+):=\lim_{t\to \alpha^+}u(t)=1, \quad  u(\beta^-):=\lim_{t\to \beta^-}u(t)=0\eeq
\beq \label{e:limsh} \lim_{t\to \alpha^+} d(u(t))|u'(t)|^{p-2}u'(t)= \lim_{t\to \beta^-} d(u(t))|u'(t)|^{p-2}u'(t)=0.\eeq
From now on, we denote by $I_u:=(\alpha,\beta)$.
%\[ I_u:=\{t\in (a,b) : 0<u(t)<1\}\]
%and we assume, without restriction, that \((a,b)\subseteq \R\) is the maximal existence interval of the solution. 
\end{defin}

Condition \eqref{e:limsh} implies that the function $(d\circ u)|u'|^{p-2}u'$ admits a continuous extension  over the entire real line as null function and allows to consider  the possible occurrence of
sharp solutions, that is solutions reaching the equilibria at
a finite time with a non-zero slope (in this case one or both the extrema $\alpha,\beta$ are
finite). However, when the existence interval is the whole real
line and $d$ is bounded in $(0,1)$, then condition \eqref{e:limsh} is automatically satisfied, as
the following result states.

\begin{prop}\label{p:1}
Let $u$ be a solution to equation in \eqref{eq:ode_lapl}, for some \(c\in \R\), satisfying \eqref{c:as}. Assume that $d$ is bounded in $(0,1)$.

Then,  if \(\alpha=-\infty\) we have  
\( \dis \lim_{t\to -\infty}d(u(t))|u'(t)|^{p-2}u'(t)=0
\). Similarly, if \(\beta=+\infty\) we have  
\( \dis \lim_{t\to +\infty}d(u(t))|u'(t)|^{p-2}u'(t)=0
\).
\end{prop}

\begin{pf} Assume that \(\beta=+\infty\).
Integrating the equation of \eqref{eq:ode_lapl} in \([0,t]\) we obtain
\begin{multline*} d(u(t))|u'(t)|^{p-2}u'(t)= 
d(u(0))|u'(0)|^{p-2}u'(0)
 - c\left(g(u(t))- g(u(0))\right) \\ + f\left(u(t)\right)- f\left(u(0)\right) - \int_0^t \rho(u(s)) {\rm d} s.
\end{multline*}
Notice that there exists the limit \(\dis\lim_{t\to +\infty} \int_0^t \rho(u(s)) {\rm d} s=\int_0^{+\infty} \rho(u(s)) {\rm d}s \in (0,
+\infty]\). So,   we infer that there exists also the limit \(\lambda:=\dis\lim_{t\to +\infty} d(u(t))|u'(t)|^{p-2}u'(t)\in [-\infty,+\infty)\). Then, since $d$ is bounded in $(0,1)$, if $\lambda\ne 0$  we get 
\( \dis \limsup_{t\to +\infty}  |u'(t)|^{p-2}u'(t) 
<0\),
implying 
that also $\dis\limsup_{t\to +\infty}u'(t)<0$, in contradiction with the boundedness of \(u\). Hence, necessarily \(\lambda=0\).
\newline
The proof concerning the limit as \(t\to -\infty\) is analogous.
\end{pf}

\md

\begin{rem}\rm 
Note that if a t.w.s. exists, then necessarily  \(c\int_0^1 g(s) \d s > \int_0^1 f(s) \d s\). Indeed, 
integrating the equation in \eqref{odeprobold} in \((a,b)\) we have
\[
\int_{a}^{b}(\dis d(u(t))|u'(t)|^{p-2} u'(t))' \d t+\int_{a}^{b}\left[cg(u(t))-f(u)\right]u'(t) \d
t+\int_{a}^{b}\rho(u(t)) \d t=0
\]
The first integral is null by \eqref{e:limsh}, so \ \[0< \int_{a}^{b}\rho(u(t)) \d t = -\int_a^b \left[cg(u(t))-f(u)\right]u'(t) \d t =  \int_0^1 \left[cg(s)-f(s)\right]  \d s.\]
Therefore,
when the dynamics does not involve convection phenomena, that if $f(u)\equiv 0$,  
then necessarily $c>0$. 
 \label{rem:1} 

\end{rem}

\md

The following result concerns the monotonicity of the solutions we
are looking for.

\md
\begin{prop}\label{p:CN}
If \ $u$ is a t.w.s. for equation \eqref{eq:RD}, then   \(u'(t)<0\) for every  \(t\in I_u\).
% and \(I_u=(\alpha,\beta)\) for some $\alpha,\beta$ with $-\infty %\le \alpha<\beta\le +\infty$.
\end{prop}

\begin{pf} 
First of all, let us prove that \(u'(t)\ne 0\) for every \(t\in I_u\). Indeed, if \(u'(t_0)=0\) for some \(t_0\in I_u\), since    \(\rho(u(t_0))>0\) we have 
\[ (d(u(t))|u'(t)|^{p-2}u'(t))'_{|t=t_0} = -\rho(u(t_0))<0\]
so the function \((d\!\circ \!u)|u'|^{p-2}u'\) is strictly decreasing in a neighborhood of \(t_0\) and it vanishes at \(t_0\).
This implies that \(u'(t)>0\) in a left neighborhood of \(t_0\).
Let 
\[t^*:=\inf \{t: u'(\tau)> 0 \text{ for every } \tau\in (t,t_0)\}.\]
By \eqref{c:as}, we have \(t^*>a\), so there exists \(u'(t^*)=0\). If \(u(t^*)>0\) we could repeat the same considerations just made for \(t_0\) and deduce that \(u'(t)<0\) in a right neighborhood of \(t^*\), in contrast to the definition of \(t^*\). So, \(u(t^*)=0\) and integrating the differential equation of \eqref{odeprobold} in \([t^*,b)\) we infer
\begin{align*} 0 & = \int_{t^*}^b (d(u(t)|u'(t)|^{p-2} u'(t))'\d t + \int_{t^*}^b [cg(u(t))-f(u(t))] u'(t) \d t + \int_{t^*}^b \rho(u(t))\d t \\ & =\int_{t^*}^b \rho(u(t))>0,
\end{align*}
a contradiction.
Therefore, we have \(u'(t)<0\) for every \(t\in I_u\).
%, implying that \(I_u\) consists of a single open interval.
\end{pf}

\bigskip

As a consequence of the monotonicity just proved, note that
  equation  \eqref{eq:ode_lapl} can be rewritten as 

\beq \label{eq:sing2}  (d(u)|u'|^{p-1})' - (c g(u)-f(u)) u' -\rho(u)=0  \eeq

\bigskip
We now recall the classification of the t.w.s.'s.
As we mentioned in Introduction, equation \eqref{eq:RD} may admit various {{types}} of solution, according to the structure of the  interval \(I_u=(\alpha,\beta)\) and the slope of \(u\) at the extreme points of \(I_u\).  More in detail, when $\alpha>-\infty$, then the t.w.s attains the equilibrium 1 at a finite time, so we can say the dynamic exhibts the phenomenon of saturation in finite times. Similarly, when $\beta<+\infty$ then the t.w.s. attains the equilibrium 0 at a finite time, hence the dynamic exhibits the phenomenon of  {\em finite speed of propagation}. Of course, the previous situations can both occur. However, notice that even when $\alpha$ and/or $\beta$ is finite the slope of the t.w.s at 1 and/or 0 could be null and in this case, the solution can be extended as a constant function, while remaining of class $C^1$.

\medskip
In order to proceed with the classification of t.w.s., recall the following definition.

\begin{defin} \rm A t.w.s. is said to be
\medskip

\begin{tabular}{ll}
{\em - \ classical} &  if \ \(\dis\lim_{t\to \alpha^+} u'(t)=\lim_{t\to \beta^-} u'(t)=0\); \\
{\em - \ sharp
of type (I)} & if \  \(\dis\lim_{t\to \alpha^+} u'(t)=0\) \ and \(\dis\lim_{t\to \beta^-} u'(t)<0\); \\
{\em - \
sharp of type (II)} & if \ \(\dis\lim_{t\to \alpha^+} u'(t)<0\) \ and \(\dis\lim_{t\to \beta^-} u'(t)=0\);\\ 
{\em - \ sharp of type (III)} & if \ \(\dis\lim_{t\to \alpha^+} u'(t)<0\) \ and \(\dis\lim_{t\to \beta^-} u'(t)<0\).
\end{tabular}

In light of the previous discussion, we point out that in the first two cases it is possible to have $\alpha>-\infty$ and in the first and third case it is possible to have  $\beta<+\infty$ (see figures in Introduction).

Clearly, if \(\dis\lim_{t\to \alpha^+} u'(t)<0\) then  \(\alpha>-\infty\) and similarly if \(\dis\lim_{t\to \beta^-} u'(t)<0\) then \(\beta<+\infty\). In this case the dynamic present the phenomenon of finite speed of propagation, since the solution reachs (or leaves) the equilibria in a finite time.
Moreover, when this happens, the solution may not be Lipschitz continuous, since the limit of the derivative can be infinite. We will see that this may happens when \(d(0)=\dot d(0)=0\) or \(d(1)=\dot d(1)=0\).

Instead, when the limits of the derivative are null, then not necessarily the corresponding extreme \(\alpha\) or \(\beta\) is infinite, but the solution admits a \(C^1-\)continuation which is constant outside the interval \((\alpha,\beta)\).

\end{defin}

\medskip

\section{An equivalent first order problem}\label{sec:ODE}

This section is devoted to introduce an equivalent singular first order problem and to recall some properties of the solutions and existence results.

\medskip
\begin{prop}\label{prop:equiv} Equation \eqref{eq:RD} admits t.w.s with speed $c$ if and only if there exists a solution to the following first order problem
\begin{equation}
\label{odeprob2}
\begin{cases}
\dot z = cg(v)-f(v)-\dfrac{(d(v))^\frac{1}{p-1}\rho(v)}{z^{\frac{1}{p-1}}} \quad \mbox{in} \;\; (0,1) \\
z(0^+)=z(1^-)=0 \\
z(v)>0 \quad \mbox{in} \;\; (0,1) 
\end{cases}
\end{equation}
{\rm(}where the notation $\dot z$ stands for derivation with respect to the variable $u${\rm )}. 

More precisely:

\begin{itemize}
\item if $u\in C^1(I_u)$  is a t.w.s. for \eqref{eq:RD}, then denoted by \(\tau=\tau(v)$ its inverse function defined in \((0,1)\),  the function $  z(v):=d(v)|u'(\tau(v))|^{p-1}$
 belongs to $C^1(0,1)$ and solves \eqref{odeprob2};

\item if $z\in C^1(0,1)$ is a solution to \eqref{odeprob2} then the (unique) solution of 
the Cauchy problem 
\beq \label{pr:cau}
\begin{cases} u'=-\left(\frac{ z(u) }{d(u)}\right)^\frac{1}{p-1} \\ u(0)=\frac12 \end{cases}\eeq
defined in its maximal existence interval \(I_u=(\alpha,\beta)\), is a t.w.s. for \eqref{eq:RD}.
\end{itemize}
\end{prop}

\begin{pf}
Let \(u\in C^1(I_u)\) be a t.w.s. for \eqref{eq:RD}. Since \(u\) is stricly decreasing in \(I_u\), with \(u'(t)<0\) for every \(t\in I_u\), there exists the inverse function \(\tau:(0,1)\to I_u\), with \(\tau\in C^1(0,1)\). Put  \(\psi_u(t):=d(u(t))|u'(t)|^{p-1}\), \(t\in I_u\),  and define
\[ z(v):= \psi_u(\tau(v)) , \quad v\in (0,1).\]
Since  \(\psi_u \in C^1(I_u)\), we deduce that \(z\in C^1(0,1)\) too. Moreover, since 
\[\psi_u'(t)= (c g(u(t))- f(u(t)) u'(t)+ \rho(u(t))\]
and
\[ u'(t)=-\left(\frac{\psi_u(t)}{d(u(t))}\right)^\frac{1}{p-1}=-\left( \frac{z(u(t))}{d(u(t))}\right)^\frac{1}{p-1},\]
for every \(v\in(0,1)\) we have
\[ \dot z(v)=  \frac{\psi_u'(\tau(v))}{u'(\tau(v))}=  (cg(v)-f(v))+\frac{\rho(v)}{u'(\tau(v))}=(cg(v)-f(v))-\frac{(d(v))^\frac{1}{p-1}\rho(v)}{(z(v))^\frac{1}{p-1}}. \]
Finally, by \eqref{e:limsh} we have
  \(\dis\lim_{v\to 0^+} z(v)=\lim_{t\to \beta^-}\psi_u(t)=0\) and   \(\dis\lim_{v\to 1^-} z(v)=\lim_{t\to \alpha^+}\psi_u(t)=0\).

\medskip
Assume now that problem \eqref{odeprob2} admits a solution \(z\in C^1(0,1)\) and consider the solution \(u\) to problem \eqref{pr:cau}, defined in its maximal existence interval \((\alpha,\beta)\).
Of course, \(u\in C^1(\alpha,\beta)\) and 
\( d(u(t))|u'(t)|^{p-1}=z(u(t))\) for every \(t\in (\alpha,\beta)\), so 
\[ ( d(u(t))|u'(t)|^{p-1})'= \dot z(u(t))u'(t) = (cg(u(t))-f(u(t))u'(t)+\rho(u(t)). \]

Moreover, since
\(u\) is strictly decreasing, there exist \(u(\alpha^+):=\dis\lim_{t\to \alpha^+} u(t)\le 1\).
If \(u(\alpha^+)<1\), then necessarily \(\alpha=-\infty\) (otherwise \((\alpha,\beta)\) is not the maximal existence interval) and 
\[0<z(u(\alpha^+)) = \lim_{t\to -\infty} z(u(t)) =\lim_{t\to -\infty} d(u(t))|u'(t)|^{p-1}=d(u(\alpha^+))\lim_{t \to -\infty} |u'(t)|^{p-1}\]
implying that there exists  \(\dis\lim_{t \to -\infty}u'(t)<0\), a contradiction. So, \(u(\alpha^+)=1\). Similarly one can prove that there exists \(u(\beta^-):=\dis\lim_{t\to \beta^-} u(t)=0\).

Finally,
\[\lim_{t\to \alpha^+}  d(u(t))|u'(t)|^{p-1}= \lim_{v\to 1^-}z(v)=0 \ \text{ and } \ \lim_{t\to \beta^-}  d(u(t))|u'(t)|^{p-1}= \lim_{v\to 0^+}z(v)=0.\] 
\end{pf}

\bigskip
A problem of the same type as \eqref{odeprob2} has been studied in \cite{Ma}, where it was considered the b.v.p. 
 \beq \begin{cases} \label{pr:sing} \dot z = cg(u)-f(u)-\dfrac{h(u)}{z^\frac{1}{p-1}} \\
z(0^+)=z(1^-)=0\\ 
z(u)>0 \ \text{ in } (0,1).
 \end{cases}\eeq
where \(h\) is continuous  in \([0,1]\), positive in \((0,1)\) and 
\(h(0)=h(1)=0\). Nevertheless, we point out that by means of a few changes, all the results proved in \cite{Ma} continue to  hold even assuming the weaker condition 
\beq h\in C[0,1)\cap L^1(0,1), \quad h(0)=0, \quad h(u)>0 \text{ in } (0,1).\label{ip:weakH}\eeq
We  now recall the main results concerning problem \eqref{pr:sing}, which will be then applied to the function  \(h(u):= (d(u))^\frac{1}{p-1} \rho(u)\), satisfying  condition \eqref{ip:weakH} as a consequence of the structural condition $(H_d)$.
The proofs are postponed in Appendix, 
where we clarify which parts of the proofs in \cite{Ma} have to be modified in light of the new weaker hypothesis \eqref{ip:weakH}.

\medskip

The first result
 deals with the maximal existence interval for the solutions to the differential equation 
\beq \label{eq:sing}
\dot z= cg(u)-f(u) -\frac{h(u)}{z^\frac{1}{p-1}}
\eeq 
  and the behavior at the extrema of the interval \((0,1)\) (see \cite[Lemma 1]{Ma}).
%From now on, with a slight abuse of notation, we denote by 
%$\dot z(u_0)$ the limit of the incremental ratio at $u_0$,  even if it is not finite. 

\bigskip

\begin{lem} 
\label{l:general} Let $h\in C[0,1)\cap L^1(0,1)$ be such that $h(u)>0$ in $(0,1)$     and let 
 \(z\in C^1(u_1,u_2)\) be a positive solution to  \eqref{eq:sing}, where \((u_1,u_2)\) is its maximal existence interval. 

Then, \(u_1=0\) and both the limits \(z(0^+)\) and \(z(u_2^-)\) exist and are finite.

\smallskip 

Moreover,  if $z(0^+)=0$ and
\beq \text{ there exists the limit } \  h_0:= \lim_{u\to 0^+}\dfrac{h(u)}{u^\frac{1}{p-1}} 
%\in [0,+\infty]
 \label{ip:h0}
\eeq 
 then  $h_0<+\infty$, $z$ is differentiable at 0 and 
  \(\dot z(0)\) is a zero of the function  
  \beq \label{def:eta0} \eta_0(t):=
t^{\frac{p}{p-1}} -(cg(0)-f(0))t^\frac{1}{p-1}+h_0, \quad t\ge 0.
\eeq

Instead, if \
\(u_2=1\), $z(1^-)=0$ and
\beq \label{ip:h1}  \text{ there exists the limit } \ h_1:= \lim_{u\to 1^-}\frac{h(u)}{(1-u)^\frac{1}{p-1}}
%\in [0,+\infty]
,\eeq
with $h_1<+\infty$, then $z$ is differentiable at 1 and 
  \(|\dot z(1)|\) is a zero of the function  
\beq \label{def:eta1}
\eta_1(t):=
t^{\frac{p}{p-1}} +(cg(1)-f(1))t^\frac{1}{p-1}-h_{1}, \quad  t\ge 0.\eeq

\end{lem}

\begin{rem}\label{r:nonex}\rm
Note the difference between the statements related to the limits $h_0$ and $h_1$: the existence of a solution to problem \eqref{pr:sing} for some $c\in \R$ is compatible with the circumstance that $h_1=+\infty$, whereas the same does not hold for $h_0=+\infty$.
%In fact, put $M:=\dis\max_{u\in [0,1]}\{|cg(u)-f(u)|:  u \in [0,1]\}$, from the equation in \eqref{pr:sing} we have 
%$\dot z(u) \le M$, so $z(u)\le Mu$ for every $u\in [0,1]$. Therefore,
%$\dot z(u)\le M- \frac{1}{M^\frac{1}{p-1}}\frac{h(u)}{u^\frac{1}{p-1}}$ for every $u\in (0,1]$. So, if $h_0=+\infty$ then $\dot z(u)\to -\infty$ as $u\to 0^+$, a contradiction.
Thus, when \eqref{pr:sing} is solvable, then necessarily $h(0)=0$ but $h$ can be unbounded near 1, provided that $h\in L^1(0,1)$.
 
Moreover, the differentiability of the solutions at 0 is always ensured, while the differentiability at 1 is guaranteed when $h_1 < +\infty$.
\end{rem}

\bigskip
In what follows we will make use of some classical results concerning differential inequalities. 
Recall that a function \(z\in C^1(I)\), with  \(I\subset (0,1)\), is  a {\em lower-solution} [resp. {\em upper-solution}] to \eqref{eq:sing} if 
\[
\dot z \le \text{\rm{[$\ge$]}} \ cg(u) -f(u) - \frac{h(u)}{z^\frac{1}{p-1}} \quad  \text{ for all } u\in I.\]

\medskip

We recall now a classical result about lower and upper-solutions (see e.g. \cite[Theorems 9.5 - 9.6]{SZ}

\begin{lem}
\label{lem:confronto}
Let \(z\) be a solution to \eqref{eq:sing}
in  \(I\subset (0,1)\)
and let \(u_0\in I\) be fixed.

If \(y\) is a lower-solution to  \eqref{eq:sing} in $(0,u_0)$ with \(z(u_0)\le y(u_0)\),  then \(z(u)\le y(u)\) for all \(u\in I\cap (0,u_0)\).

If \(y\) is an upper-solution for the equation of \eqref{eq:sing} in $(0,u_0)$ with \(z(u_0)\ge y(u_0)\),  then \(z(u)\ge y(u)\) for all \(u\in I\cap (0,u_0)\).

\end{lem}

The following comparison criterium will be  useful  in order to establish  the existence of solutions to problem \eqref{pr:sing}  (see \cite[Proposition 1]{Ma}).

\begin{prop} Assume that \eqref{ip:weakH} holds true. 
Suppose there exists a lower-solution \(\varphi\) to  \eqref{eq:sing} in the whole interval \((0,1)\), such that  
 \(\varphi(0^+)=0\) and \(\varphi(u)>0\) for every \(u\in (0,1)\).
 
Then, the singular boundary value problem \eqref{pr:sing} admits a  solution \(z\), such that \(0<z(u)<\varphi(u)\) for every \(u\in (0,1)\).
 \label{l:soprasol}
\end{prop}

Finally, we now report the existence theorem (see \cite[Theorem 1]{Ma}).

\begin{thm}\label{t:main} Let  \eqref{ip:h0} be satisfied, with $h_0<+\infty$.

Then, 
there exists a threshold value \(c^*\) such that problem \eqref{pr:sing} admits a solution if and only if \(c\ge c^*\).
Moreover, put 
\[ G_0:= \inf_{u\in (0,1)} \fint_0^u g(s) \ \d s\ , \quad F_0:=\sup_{u\in (0,1)} \fint_0^u f(s) \ \d s  \ , \quad  
H_0:= \sup_{u\in (0,1)} \fint_0^u\frac{h(s)}{s^\frac{1}{p-1}} \d s,
\]
(where \(\fint\) stands for the mean  value) we have

\[
\frac{f(0)}{g(0)} + \frac{p'(p-1)^\frac{1}{p}}{g(0)}\ h_0^\frac{1}{p'}\le c^*\le \frac{F_0}{G_0}+ \frac{ p'(p-1)^\frac{1}{p}}{G_0}\ H_0^\frac{1}{p'}.
\]
Finally, for every \(c\ge c^*\) the solution is unique.

\end{thm}

\begin{rem} \label{r:eta0}
\rm
Observe that when a solution exists for some $c\in \R$ and \eqref{ip:h0} is fulfilled,  then by Lemma \ref{l:general} we deduce that $\dis\min_{t\ge 0} \eta_0(t)\le 0$ (see \eqref{def:eta0}). Moreover, simple calculations show that $\dis \min_{t\ge 0} \eta_0(t)<0$ if and only if 
\[ cg(0)-f(0) > 
 p'(p-1)^\frac{1}{p}\ h_0^\frac{1}{p'}.\]
 Finally, 
it is immediate to verify that when $\dis \min_{t\ge 0} \eta_0(t) < 0$  there exist \(r_0^-, r_0^+\), with  \(0\le r_0^-< r_0^+\), such that \(\eta_0(t) <0 \) if and only if \( r_0^- < t < r_0^+\).

\end{rem}

We conclude this section devoted to the first order boundary value problem \eqref{pr:sing} with two results concerning the slope of the solutions at the extremes of the interval $(0,1)$.

\medskip

\begin{lem}
\label{lem:zderiv0} Let \eqref{ip:h0} be satisfied, with $h_0<+\infty$. Assume that problem \eqref{pr:sing} admits a solution \(z_c\) for some \(c\in \R\).
Let   \(\eta_0\) be the function defined in \eqref{def:eta0}.
 
 Then 
\[ \dot z_c(0)=\begin{cases}r_0^- & \text{ if } c>c^* \\ 
r_0^+ & \text{ if } c=c^* \end{cases} \]
where $r_0^-\le r_0^+$ are the zeros (possible coincident) of the function $\eta_0$ {\rm (}see Remark \ref{r:eta0}{\rm )}.
 
\end{lem}  

\begin{pf} 
Since problem \eqref{pr:sing} admits a solution for the given \(c\), by  Lemma \ref{l:general}  we have that \(z_c\) is differentiable at \(0\), and $\eta_0(\dot z_c(0))=0$. So, $\min \eta_0(t) \le 0$ and 
\( \dot z_c(0)\in \{r_0^-, r_0^+\}\).

\medskip 
Let us first consider the case \(c>c^*\) and assume, by contradiction, 
\(\dot z_c(0)=r_0^+\).
Since \(g(0)>0\), the function \(c\mapsto r_0^{+}(c)\) is strictly increasing, then  \( r_0^-(c^*)\le r_0^+(c^*) < r_0^+(c) \). Therefore,   
\(\dot z_c(0) > \dot z_{c^*}(0)\), implying \(z_c(u)>z_{c^*}(u)\) in a right neighbourhood of \(0\).
Put
\[ u_0:=\sup \{ u\in (0, 1): \ z_c(\xi) > z_{c^*}(\xi) \text{ for every } \xi\in (0,u) \}.\]
Since \(z_c(1^-)=z_{c^*}(1^-)=0\), we have \(z_c(u_0^-)=z_{c^*}(u_0^-)\) and  by \eqref{ip:g} and \eqref{pr:sing} we deduce
\[ 0= z_c(u_0^-)-z_{c^*}(u_0^-)=(c-c^*)\int_0^{u_0}\!\! g(s)  \d s - \int_0^{u_0}\!\! h(s)\left(\frac{1}{(z_c(s))^\frac{1}{p-1}} - \frac{1}{(z_{c^*}(s))^\frac{1}{p-1}}\right) \d s >0,\]
a contradiction. So, \(\dot z_c(0)=r_0^-\) when \(c>c^*\).

\medskip
Let us now consider the case \(c=c^*\). If \(r_0^-=r_0^+\) there is nothing to prove.  So, from now on assume that \(r_0^-<r_0^+\) and suppose, by contradiction, \(\dot z_{c^*}(0)=r_0^-\).

 Fix a value \(\theta\) with \(r_0^-<\theta<r_0^+\), so that \(\eta_0(\theta)<0\) (see Remark \ref{r:eta0}) and let \(\varepsilon>0\) be such that \( \eta_0(\theta)+\varepsilon \theta^\frac{1}{p-1} <0\).
Let us now prove that there exists a value \(\rho>0\) such that 
\beq
cg(u)-f(u) - \frac{h(u)}{(\theta u)^\frac{1}{p-1}} > \theta \quad \text{whenever } |c-c^*|<\rho, \ u<\rho.\label{e:claim}
\eeq
Indeed, 
by the continuity of the two-variables function
\((c,u)\mapsto cg(u)-f(u)\) at the point  \((c^*,0)\), there exists \(\rho>0\) such that if \(|c-c^*|<\rho\) and \(0<u<\rho\) we have
\[  cg(u)-f(u) > cg(0)-f(0) - \frac{\varepsilon}{2}.\] 
Moreover, by \eqref{ip:h0} we can assume without restriction that 
\[ \frac{h(u)}{(\theta u)^\frac{1}{p-1}} <  \frac{h_0}{\theta^\frac{1}{p-1}} +\frac\varepsilon{2}  \quad \text{whenever } u<\rho.\]
So, if \(|c-c^*|<\rho\) and \(0<u<\rho\), by the choice of $\varepsilon$ we have
\[  cg(u)-f(u)-\frac{h(u)}{(\theta u)^\frac{1}{p-1}} > cg(0)-f(0)-\frac{h_0}{\theta^\frac{1}{p-1}} - \varepsilon= \theta- \frac{1}{\theta^\frac{1}{p-1}}\eta_0(\theta) - \varepsilon>\theta, \] 
that is \eqref{e:claim}. 

\bigskip
Therefore, if we consider the function \(\lambda(u):= \theta u\), by \eqref{e:claim} we infer it is a  lower-solution to \eqref{eq:sing} in the interval \((0,\rho)\) for every \(c\) satisfying \(|c-c^*|<\rho\).
Moreover, since \(\theta>r_0^-=\dot z_{c^*}(0)\), we deduce that \(
\theta u> z_{c^*}(u)\) in \((0,\delta]\), for some \(\delta<\rho\).

\medskip
Let us now consider the solution \(\varphi^*\) to \eqref{eq:sing} for \(c=c^*\), passing throughout \((\delta, \theta \delta)\), defined in its maximal existence interval \((u_1,u_2)\subseteq (0,1)\). By Lemma \ref{l:general} we have \(u_1=0\) and there exist, finite, 
  the values \(\varphi^*(0^+)\) and \(\varphi^*(u_2^-)\). 
  
Since \(\lambda\) is a lower-solution and $\lambda(\delta)=\theta\delta=\varphi^*(\delta)$, by Lemma \ref{lem:confronto} we obtain \(0\le \varphi^*(u)\le \lambda(u)=\theta u\) for every \(u\in (0,\delta)\), so  \(\varphi^*(0^+)=0\). 
   Moreover, since \(\varphi^*(\delta)>z_{c^*}(\delta)\), by the uniqueness of the solution to  \eqref{eq:sing} passing through a given point, we have \(\varphi^*(u) > z_{c^*}(u)\) for every \(u\in [\delta,u_2)\). So, \(u_2=1\).

\medskip
Observe now that \(\varphi^*(1^-)>0\). 
Indeed, if \(\varphi^*(1^-)=0\), then
 \(\varphi^*\) is a solution to problem \eqref{pr:sing} for \(c=c^*\), different from \(z_{c^*}\), and this contradicts the uniqueness of the solution of problem \eqref{pr:sing}. 

By the continuous dependence on the data, there exists a value \(\tilde c \in (c^*-\rho,c^*)\) such that the solution \(\tilde \varphi\) to  \eqref{eq:sing} for \(c=\tilde c\), passing through \((\delta,\theta \delta)\) is defined in \([\delta,1]\), with \(\tilde \varphi(1)>0\). 
Moreover, since \(\lambda\) is a lower solution in \((0,\delta]\) also for \(c=\tilde c\), we deduce that \(0\le \tilde\varphi(u)\le \lambda(u)=\theta u\) for every \(u\in (0,\delta]\), so \(\tilde\varphi(0^+)=0\). Therefore, 
 \(\tilde\varphi\) satisfies the assumptions of Proposition \ref{l:soprasol} and then problem \eqref{pr:sing} admits a solution even for \(c=\tilde c\), in contrast to the definition of \(c^*\) as the minimum admissible wave speed.
\end{pf}

\bigskip

In order to prove an analogous result regarding the slope of the solutions at 1, we will use the following comparison result.

\begin{lem} \label{l:comp1} Let \eqref{ip:weakH} be satisfied. 
Let $z$ be a solution of \eqref{pr:sing} and let $y\in C^1(1-\delta,1)$ be a positive lower-solution  {\rm [}resp. upper-solution{\rm]} to \eqref{eq:sing}, with $y(1^-)=0$. 

Then, $y(u)\ge z(u)$ {\rm[}resp. $y(u)\le z(u)${\rm]} for every $u\in (1-\delta,1)$.

\end{lem}

\begin{pf}
We provide the proof in the case $y$ is a lower-solution (the proof in the other case is analogous).

Assume, by contradiction, $y(u_0)<z(u_0)$ for some $u_0\in (1-\delta,1)$. Let us define $u_1:=\sup\{ u\in (u_0,1): \ y(s)<z(s) \text{ for every } s\in [u_0,u)\}$.
Then, 
\begin{align*} z(u_1^-)= & z(u_0)+\int_{u_0}^{u_1} \left(cg(u)-f(u)- \frac{h(u)}{(z(u))^\frac{1}{p-1}} \right)\d u \\ >& y(u_0) + \int_{u_0}^{u_1} \left(cg(u)-f(u)- \frac{h(u)}{(y(u))^\frac{1}{p-1}}\right) \d u \\ \ge & y(u_0)+ \int_{u_0}^{u_1} \dot y(u) \d u = y(u_1^-)  \end{align*}
yielding $u_1=1$ and consequently $z(1^-)>0$, a contradiction.  
\end{pf}

\bigskip

Now we discuss the behavior of the solutions at 1. As we have already observed in Remark \ref{r:nonex}, while differentiability in 0 is always ensured, the same is not true for differentiability in 1. 
%So, from now on, with a slight abuse of notation, we denote by $\dot z(1)$ the limit of the incremental ratio at $u=1$, even it is not finite (provided that it exists).

\begin{lem} \label{lem:zderiv1}
 Assume that problem \eqref{pr:sing} admits a (unique) solution \(z_c\) for some \(c\in \R\) and let \eqref{ip:h1} be satisfied. Then, 
 $z_c$ is differentiable at $1$ if and only if $h_1<+\infty$.

Moreover, if \(0<h_1<+\infty\) then \(\dot z_c(1)<0\); while if \(h_1=0\) then we have $\dot z_c(1)=\min \{0, cg(1)-f(1)\}$.
\end{lem}

\begin{pf} First assume $h_1=+\infty$. If $z$ is differentiable at 1,   from \eqref{eq:sing} we deduce $\dot z(u)\to -\infty$ as $u\to 1^-$, a contradiction. Vice versa, if $h_1<+\infty$ then  $z$ is  differentiable at 1, as a consequence of Lemma \ref{l:general}.

Moreover,  when $h_1<+\infty$, we can consider the function
\(\eta_1(t)\) introduced in \eqref{def:eta1}.  Of course, $\eta_1(0)=-h_1$, so,  when \(h_1>0\)  we have \(\dot z(1)<0\)  by Lemma \ref{l:general}.

Assume now \(h_1=0\). In this case, if $cg(1)-f(1)\ge 0$, the function \(\eta_1\) vanishes only at \(t=0\), so by Lemma \ref{l:general} we get \(\dot z(1)=0\). Instead,  if  $cg(1)-f(1)< 0$, then \(\eta_1\)  vanishes at 0 and \(f(1)-cg(1)\).  But for \(u\) sufficiently close to 1 we have
\[ \dot z(u)=cg(u)-f(u) - \frac{h(u)}{(z(u))^\frac{1}{p-1}} < cg(u)-f(u) < \frac12 (cg(1)-f(1))<0\]
hence necessarily $\dot z(1)= cg(1)-f(1)$. 

Summarizing, when \(h_1=0\) then  $\dot z_c(1)=\min \{0, cg(1)-f(1)\}$.

\end{pf}

\section{Existence and uniqueness result}

As a consequence of Proposition \ref{prop:equiv}, Lemma \ref{l:general} and Theorem \ref{t:main}, we can now state the existence result for solutions of problem \eqref{odeprobold}.

\begin{thm}\label{t:main2} 
Assume that  
\beq \text{ there exists the limit } \  \ell_0:= \lim_{u\to 0^+}\rho(u)\dfrac{(d(u))^\frac{1}{p-1}}{u^\frac{1}{p-1}} \in [0,+\infty]. \label{ip:l0}
\eeq

Then, if \(\ell_{0}=+\infty\) problem \eqref{odeprobold} does not admit solutions for any \(c\in \R\).
Otherwise, if \(\ell_{0}<+\infty\)
there exists a threshold value \(c^*\) such that problem \eqref{pr:sing} admits solution if and only if \(c\ge c^*\).
Moreover, put 
\[ G_0:= \inf_{u\in (0,1)} \fint_0^u g(s)  \d s\ , \quad F_0:=\sup_{u\in (0,1)} \fint_0^u f(s)  \d s  \ , \quad  
L_0:= \sup_{u\in (0,1)} \fint_0^u\left(\frac{d(s)}{s}\right)^\frac{1}{p-1}\!\!\!\!\rho(s) \d s,
\]
(where \(\fint\) stands for the mean  value) we have

\beq
\frac{f(0)}{g(0)} + \frac{p'(p-1)^\frac{1}{p}}{g(0)}\ \ell_0^\frac{1}{p'}\le c^*\le \frac{F_0}{G_0}+ \frac{ p'(p-1)^\frac{1}{p}}{G_0}\ L_0^\frac{1}{p'}.
\label{eq:stimanuova}
\eeq
Finally, for every \(c\ge c^*\) the solution is unique.

\end{thm}

It is worth noting that estimate \eqref{eq:stimanuova} 
reduces to that already obtained in \cite{mp1} for the classical case ($p=2$ and $g(u)\equiv 1$). Moreover, 
notice that \eqref{eq:stimanuova} implies $c^*g(0)-f(0)\ge 0$. When $f(u)\equiv 0$, by  Remark \ref{rem:1} and \eqref{ip:g} we have $c^*>0$, but if $f$ is not null it may occur that $c^*g(0)=f(0)$.  As we will show in the next section, in order to classify the t.w.s. with respect to their behavior at the equilibrium 0, it is important to establish if $c^*g(0)-f(0)> 0$ or $c^*g(0)-f(0)= 0$. So, this issue deserves a further study. 
The next result provide sufficient criteria ensuring the strict inequality or the equality. This theorem takes up the idea underlying the analogous result proved in \cite[Theorems 1.2, 1.3]{MM1} for the case $p=2$ and $g(u)\equiv 1$.

\begin{thm} 
\label{t:stima}
Let \eqref{ip:l0} be fulfilled, with $\ell_0<+\infty$, and  let $k\in \R$ be fixed. 

The following statements hold true:

\smallskip
{\bf \rm (a)} \  if  $f(u)\ge kg(u) $ in a right neighborhood of 0, then $c^*>k$;

\smallskip
{\bf \rm (b)} \ if for every $u\in (0,1)$ we have 
 $kg(u)\ge f(u)$  and 
\beq
\label{ip:cond0}
\left(kg(u)-f(u) \right)^{p-1}  \int_0^u (kg(s)-f(s))\d s \ge \frac{p^p}{(p-1)^{p-1}} d(u)(\rho(u))^{p-1}
\eeq
then, $c^* \le k$;

\smallskip
{\bf \rm (c)} \ if there exists  $\delta>0$ such that in $(0,\delta)$ we have    $ kg(u)\ge f(u)  $  and 
\beq\label{ip:cond0bis} \left(kg(u)-f(u) \right)^{p-1}  \int_0^u (kg(s)-f(s))\d s \le \ell  d(u)(\rho(u))^{p-1}\eeq
 for some $\ell < \frac{p^p}{(p-1)^{p-1}}$, then $c^*>k$.

\end{thm}

\begin{pf}
The assertion $(a)$ is trivial, since if $c^*\le k$ then the solution $z$ of problem \eqref{odeprob2} for $c=k$ satisfies
\[ \dot z(u) = kg(u)- f(u) - \frac{\rho(u)(d(u))^\frac{1}{p-1}}{(z(u))^\frac{1}{p-1}}\le 0\]
in a right neighborhood of 0, a contradiction since $z(0^+)=0$ and $z(u)>0$ in $(0,1)$.

As for assertion $(b)$, put 
\[ \varphi(u):=  \frac{1}{p}(kg(u)-f(u)), \quad \Phi(u):= \int_0^u 
\varphi(s) \d s .  \]
Then,  by \eqref{ip:cond0} for every $u\in (0,1)$ we have  $\Phi(u)>0$ and 
\[
\frac{\rho(u)d(u)^\frac{1}{p-1}}{(\Phi(u))^\frac{1}{p-1}}= p^\frac{1}{p-1} \frac{\rho(u)d(u)^\frac{1}{p-1}}{\left(\int_0^u 
(kg(s)-f(s)) \d s \right)^\frac{1}{p-1}} \le\]
\[\le \frac{p-1}{p}(kg(u)-f(u)) = kg(u) - f(u) -  \dot \Phi(u), \]
so $\Phi$ is a positive lower-solution to \eqref{eq:sing} for $c=k$. Therefore, by Proposition \ref{l:soprasol} problem \eqref{odeprob2} admits a solution for $c=k$. So,   $c^*\le k$.

\medskip
Finally, let us prove assertion $(c)$.
 Suppose by contradiction that problem \eqref{odeprob2} admits a solution  $z$ for $c=k$.

Consider again the functions $\varphi$ and $\Phi$ above defined, which are positive in $(0,\delta)$. Then, by \eqref{ip:cond0bis}, in $(0,\delta)$ we have 
\beq 
\frac{\rho(u)(d(u))^\frac{1}{p-1}}{(\Phi(u))^\frac{1}{p-1}}\ge \frac{p^\frac{p}{p-1}}{ \ell^\frac{1}{p-1}} \varphi(u). \label{eq:eq4}
\eeq

Observe now that if $\dot z(u)\le \alpha p  \varphi(u)$ in $(0,\delta)$ for some $\alpha>0$, then  $z(u)\le \alpha p \Phi(u)$, implying 
\begin{align*} \dot z(u) = p\varphi(u) - \frac{ (d(u))^\frac{1}{p-1} \rho(u)}{(z(u))^\frac{1}{p-1}} & \le p\varphi(u) - \frac{ (d(u))^\frac{1}{p-1} \rho(u)}{(\alpha p \Phi(u))^\frac{1}{p-1}} \\ & \stackrel{\eqref{eq:eq4} }{\le} p\varphi(u) - \frac{p}{(\ell\alpha )^\frac{1}{p-1}}\varphi(u)  =\left(1-(\ell\alpha)^{-\frac{1}{p-1}}\right) p\varphi(u). \end{align*}

\noindent
So, let us consider now the sequence defined by $a_1:=1$, $a_{n+1}:=1- \left(\ell a_n \right)^{-\frac{1}{p-1}}$, until $a_n>0$. By what we have just proved, since $\dot z(u)\le a_1p \varphi(u)$ in $(0,\delta)$, we have $\dot z(u)\le a_{n+1} p \varphi(u) $ in $(0,\delta)$, for every $n$ such that  $a_n>0$.

Put $\mu(t):= t-1+(\ell t)^{-\frac{1}{p-1}}$, 
since $\ell < \frac{p^p}{(p-1)^{p-1}}$ one can verify that 
$\mu(t)>0$ for every $t>0$.
 Then, as long as $a_n>0$, we have $a_{n+1}<a_n$.
Hence, if $a_n>0$ for every $n\in \N$, we should derive $a_n>1/\ell$ for every $n\in \N$ and $a_n\to a^*$, for some $a^*>0$ with $\mu(a^*)=0$, a contradiction. 
 Therefore, we have $a_{\bar n}>0$ and $a_{\bar n+1}<0$ for some $\bar n \in \N$, implying  $\dot z(u)<
a_{\bar n+1}p \varphi(u)<0$ in $(0,\delta)$, 
 a contradiction.
\end{pf}

\bigskip

\begin{rem}\label{rem:c*}\rm
Putting $k:=f(0)/g(0)$ in statement $(b)$ of Theorem \ref{t:stima} we derive a sufficient condition to have $c^*g(0)=f(0)$, since from \eqref{eq:stimanuova}  it follows $c^*g(0)\ge 0$. Similarly,  statement $(c)$ provides a sufficient condition to have  $c^*g(0)>f(0)$.

\end{rem}

\section{Local analysis at the equilibrium 0}

In this section we study the property of {\em finite speed of propagation} of equation \eqref{eq:RD}, occurring when the solution reaches the equilibrium 0 in a finite time, that is the extremum $\beta$ of the existence interval  is finite.

\medskip
First  of all, observe that 
\beq \label{eq:integr}
\beta= -\int_0^{\beta}\!\!\left(\frac{d(u(t))}{z(u(t))} \right)^\frac{1}{p-1} \!\! u'(t)\d t= \int_0^\frac12 \left(\frac{d(u)}{z(u)}\right)^\frac{1}{p-1}\!\!  \d u,\eeq
so, \(\beta\) is finite or infinite according to the convergence or the divergence of the last integral in \eqref{eq:integr}.

\medskip
The first result concerns the t.w.s. having speed $c>c^*$. We show that in this case the convergence/divergence of the last integral in \eqref{eq:integr} is uniquely related to the convergence/divergence of the integral \(\dis\int_0^\frac12 \frac{1}{\rho(u)} \d u\).

\bigskip
\begin{thm}
\label{t:T2} 
Let \eqref{ip:l0} be fulfilled, with $\ell_0<+\infty$.
Let \(u\) be a solution of \eqref{odeprobold} for some \(c> c^*\).
Assume that 
%\eqref{ip:chain} holds true.
 there exist three constants \(\delta, \mu, M>0\) and a  function $\xi\in C^1[0,\delta]$  such that 
\beq \label{ip:chain}
\mu\xi(u) \le  d(u)(\rho(u))^{p-1}\le M\xi(u) \quad \text{ for every } u\in (0,\delta).
\eeq

Then, $\beta=+\infty$ if and only if $\dis\int_0^\frac12 \frac{1}{\rho(u)} \d u=+\infty$. Morever, even if $\beta<+\infty$ we have $ \dis\lim_{t\to \beta^-} u'(t)=0$.
\end{thm}

\noindent
\begin{pf} Let $z$ be the  solution of \eqref{odeprob2} corresponding to the t.w.s. $u$ (see Proposition \ref{prop:equiv}).
Recall that by Lemma \ref{l:general} $z$ is differentiable at 0, with $\eta_0(\dot z(0))=0$  (see \eqref{def:eta0}).

\medskip
\(\bullet\) {\bf Case 1:} \(\ell_0>0\).

\smallskip

In this case, since $\eta_0(0)>0$,  we have $\dot z(0)>0$. So, as $u\to 0$ we have
\[ \left( \frac{d(u)}{z(u)}\right)^\frac{1}{p-1} \sim \left( \frac{d(u)}{\dot z(0) u}\right)^\frac{1}{p-1} \sim \left(\frac{\ell_0}{\dot z(0)}\right)^\frac{1}{p-1}\!\! \frac{1}{\rho(u)}\]
yielding the assertion by \eqref{eq:integr}.

\medskip
\(\bullet\) {\bf Case 2:} \(\ell_0=0\).
\smallskip

In this case, since $c>c^*$ and $\eta_0(0)=0$, by Lemma \ref{lem:zderiv0} we have  \(\dot z(0)=r_0^-=0\). Then, there exists a sequence \((u_n)_n\) converging to 0, such that \(\dot z(u_n)\to 0\). Hence, put $\lambda_0:=cg(0)-f(0)>c^*g(0)-f(0)\ge 0$,
for \(n\) sufficiently large we have 
\beq  \frac12 \lambda_0 < 
 cg(u_n)-f(u_n)-\dot z(u_n) < \frac32 \lambda_0. \label{eq:sequence}
\eeq
 Put 
\[k_0:=\left(\frac{2}{\lambda_0}\right)^{p-1} \ , \quad
\psi(u):=Mk_0\xi(u).\]
Of course, \(\psi\in C^1([0,\delta])\), with 
$\dot \psi(0)=0$ by \eqref{ip:chain}, since $\ell_0=0$. Moreover, again by \eqref{ip:chain} we have that 
\[ \limsup_{u \to 0} \rho(u) \!\left(\frac{d(u)}{\psi(u)}\right)^\frac{1}{p-1}\le \frac{1}{k_0^\frac{1}{p-1}}=\frac12 \lambda_0.\]  
Therefore, we have
\begin{align*}  \dot\psi(0) & =0 < \frac12\lambda_0 \le
\liminf_{u\to 0} \left(cg(u)-f(u) - \frac{(d(u))^\frac{1}{p-1} \rho(u)}{(\psi(u))^\frac{1}{p-1}}\right).\end{align*}
Then, for some \(\rho\in(0,\delta)\) we have 
\[ \dot\psi(u) < cg(u)-f(u) - \frac{(d(u))^\frac{1}{p-1} \rho(u)}{(\psi(u))^\frac{1}{p-1}} \quad \text{ for every } u \in (0,\rho)\]
i.e. \(\psi\) is a lower-solution in the interval \((0,\rho)\).

\smallskip
On the other hand, for   \(n\) large enough  we have
\begin{align*} \left(\frac{\psi(u_n)}{z(u_n))}\right)^\frac{1}{p-1} & = (Mk_0)^\frac{1}{p-1}\left(\frac{\xi(u_n)}{z(u_n)}\right)^\frac{1}{p-1}\stackrel{\eqref{ip:chain}}{ \ge} k_0^\frac{1}{p-1}\rho(u_n)\left(\frac{d(u_n)}{z(u_n)}\right)^\frac{1}{p-1}\\ &= k_0^\frac{1}{p-1}(cg(u_n)-f(u_n)-\dot z(u_n))\stackrel{\eqref{eq:sequence}}{>}1. \end{align*}
Hence,  there exists a value \(\varepsilon<\rho\) such that 
 \(z(\varepsilon)<\psi(\varepsilon)\). So, being \(\psi\) a lower-solution in $(0,\varepsilon]$, by Lemma \ref{lem:confronto} we infer \(z(u)\le \psi(u)\) for every \(u\in (0,\varepsilon)\). Thus, we 
   deduce that  
\beq \left(\frac{d(u)}{z(u)} \right)^\frac{1}{p-1}\ge  \left(\frac{d(u)}{\psi(u)} \right)^\frac{1}{p-1} \stackrel{\eqref{ip:chain}}{\ge}\left(\frac{\mu}{Mk_0}\right)^\frac{1}{p-1}\!\!\frac{1}{\rho(u)} \quad \text{ for every } u\in (0,\varepsilon). \label{eq:stimasotto} \eeq

\medskip

Put now  $\varphi(u):=\dfrac{1}{3^{p-1}}\mu k_0\xi(u)$.
Notice that 
%$\psi(\delta)<z(\delta)$ and
\[\dot \varphi(0)=0>\lambda_0-\frac{3}{k_0^\frac{1}{p-1}}\stackrel{\eqref{ip:chain}}{\ge}
\limsup_{u \to 0} \left(cg(u)-f(u)- \frac{(d(u))^\frac{1}{p-1}\rho(u)}{(\varphi(u))^\frac{1}{p-1}}\right).
\]
Therefore, for some $\sigma\in (0,\delta)$ we have
\[\dot \varphi(u)> cg(u)-f(u)- \frac{(d(u))^\frac{1}{p-1}\rho(u)}{(\varphi(u))^\frac{1}{p-1}} \quad \text{ for every } u\in (0,\sigma),\]
that is $\varphi$ is an upper-solution in $(0,\sigma)$.

On the other hand,  for $n$ large enough we have
%\[cg(u_n)-f(u_n)-\dot z_c(u_n)<\frac{3}{k_0^\frac{1}{p-1}}\]
%and thus
\begin{align*} \left(\frac{\varphi(u_n)}{z(u_n))}\right)^\frac{1}{p-1} & =\frac13(\mu k_0)^\frac{1}{p-1}\left(\frac{\xi(u_n)}{z(u_n)}\right)^\frac{1}{p-1} \stackrel{\eqref{ip:chain}}{ \le}\frac13k_0^\frac{1}{p-1}\rho(u_n)\left(\frac{d(u_n)}{z(u_n)}\right)^\frac{1}{p-1}\\ &= \frac13 k_0^\frac{1}{p-1}(cg(u_n)-f(u_n)-\dot z(u_n))\stackrel{\eqref{eq:sequence}}{<}1. \end{align*}
Therefore, we have $ z(\eta) >\varphi(\eta)$ for some $\eta\in (0,\sigma)$, implying that $z(u)\ge \varphi(u)$ in $(0,\eta]$, since $\varphi$ is an upper-solution (see Lemma \ref{lem:confronto}).
Thus, 
\beq\left(\frac{d(u)}{z(u)} \right)^\frac{1}{p-1}\le 
\frac{3}{(\mu k_0)^\frac{1}{p-1}}
 \left(\frac{d(u)}{\xi(u)} \right)^\frac{1}{p-1} \stackrel{\eqref{ip:chain}}{\le}3
 \left(\frac{M}{\mu k_0}\right)^\frac{1}{p-1}\!\!
 \frac{1}{\rho(u)} \quad \text{in } (0, \eta). \label{eq:stimasopra} \eeq
Therefore, the assertion converning the boundedness  or unboundedness of $\beta$  follows from \eqref{eq:stimasotto} and \eqref{eq:stimasopra}.
Finally, note that by \eqref{eq:stimasotto} we have 
\(  \left( \frac{z(u)}{d(u)}\right)^\frac{1}{p-1} \le  \rho(u) \left(\frac{Mk_0}{\mu} \right)^\frac{1}{p-1} \rho(u) \), so $\dis\lim_{t\to \beta^-}u'(t) =0$ by \eqref{pr:cau}.
\end{pf}

\begin{rem}
\rm
In view of the proof of the previous theorem, observe that assumption \eqref{ip:chain} can be removed when $\ell_0>0$.

 Moreover, in the case $\ell_0=0$, if  the function $t\mapsto d(u)\rho(u)^{p-1}$  is $C^1$ in an open right neighborhood of 0 and there exists (finite or infinite) the limit  of its derivative as $u\to 0$,
 then it is necessarily null and  assumption \eqref{ip:chain} is trivially satisfied taking $\xi(u):=d(u)(\rho(u))^{p-1}$. 

Finally, note that \eqref{ip:chain} is fulfilled even if  $d(u)(\rho(u))^{p-1}\sim \xi(u)$ as $u\to 0$, for some $C^1$-function $\xi$, defined in a right neighborhood $[0,\delta)$. For instance, \eqref{ip:chain} holds true when  \(d(u)(\rho(u))^{p-1}\sim Cu^s$ for some $C,s>0$. Then, Theorem \ref{t:T2} provides the answer to the open problems in \cite{DZ} (see also \cite{DJKZ}) where the classification result has been obtained only for sufficiently large $c>c^*$.

\end{rem}

\medskip
The following result provides the   classification of the t.w.s. having speed $c^*$.

\begin{thm}
\label{t:T2*} 
Let \eqref{ip:l0} be fulfilled, with $\ell_0<+\infty$.
Let \(u\) be a solution to \eqref{odeprobold} for \(c= c^*\).

Then:

\smallskip
{\bf \rm (a)} \
 if \ $\ell_0>0$, then
 \[ \beta<+\infty \quad \text{ if and only if } \quad  \dis \int_0^\frac12 \frac{1}{\rho(u)}\d u<+\infty.\]  
 
\smallskip
{\bf \rm (b)} \
 if \ $\ell_0=0$ and  $c^*g(0)>f(0)$,  we have 
\[ \beta<+\infty\quad \text{ if and only if } \quad \dis\int_0^\frac12 \!\!\left( \frac{d(u)}{u}\right)^\frac{1}{p-1}\!\! \d u<+\infty;\]

\smallskip
{\bf \rm (c)} \
 if $\ell_0= c^*g(0)-f(0) =0$ and \eqref{ip:chain} holds true, then
 \[ \dis \int_0^\frac12 \frac{1}{\rho(u)}\d u<+\infty \quad \text{ implies } \quad \beta<+\infty.\]  
\end{thm}

\begin{pf}
 Let $z$ be the  solution of problem \eqref{odeprob2} corresponding to the t.w.s. $u$ (see Proposition \ref{prop:equiv}).
By virtue of Lemma \ref{lem:zderiv0},  we have $\dot z(0)=r_0^+$, where $r_0^+$ is the biggest zero of the function $\eta_0$ (see \eqref{def:eta0}).
 Then,  if $\ell_0>0$  or  $c^*g(0)>f(0)$, we have $\dot z_*(0)>0$.
Therefore, if $\ell_0>0$ we have
\[  \frac{\ell_0}{\rho(u)} \sim \left(\frac{d(u)}{u}\right)^\frac{1}{p-1}\sim  \left(\frac{\dot z(0) d(u)}{ z(u)}\right)^\frac{1}{p-1} \quad \text{as } u \to 0\]
and the assertion follows from \eqref{eq:integr}.

Instead, if $\ell_0=0$ and  $c^*g(0)>f(0)$ we have
	\(  \dfrac{d(u)}{z(u)}\sim  \dfrac{d(u)}{\dot z(0)u} \quad \text{as } u \to 0\),
and the assertion again follows from \eqref{eq:integr}.

Finally, if  $\ell_0= c^*g(0)-f(0) =0$, by Lemma  \ref{lem:zderiv0} we have $\dot z(0)=0$, since $\eta_0(t)>0$ for every $t>0$.
 Then, there exists a sequence \((u_n)_n\) converging to 0, such that \(\dot z(u_n)\to 0\). Hence, 
for \(n\) sufficiently large we have 
\beq  
 c^*g(u_n)-f(u_n)-\dot z(u_n) < \mu^\frac{1}{p-1}, \label{eq:seq*}
\eeq
where $\mu$ is the constant involved in \eqref{ip:chain}.

Moreover, 
$\dot \xi(0)=0$ by \eqref{ip:chain}, since $\ell_0=0$. Thus, again by \eqref{ip:chain} we have 
\begin{align*}  \dot\xi(0) & =0 > -\mu^\frac{1}{p-1} \ge
\limsup_{u\to 0} \left(c^*g(u)-f(u) - \frac{(d(u))^\frac{1}{p-1} \rho(u)}{(\xi(u))^\frac{1}{p-1}}\right).\end{align*}
Then, for some \(\rho\in(0,\delta)\) we have 
\[ \dot\xi(u) > c^*g(u)-f(u) - \frac{(d(u))^\frac{1}{p-1} \rho(u)}{(\xi(u))^\frac{1}{p-1}} \quad \text{ for every } u \in (0,\rho)\]
i.e. \(\xi\) is an  upper-solution in the interval \((0,\rho)\).

\smallskip
On the other hand, for   \(n\) large enough  we have
\begin{align*} \left(\frac{\xi(u_n)}{z(u_n))}\right)^\frac{1}{p-1}  \stackrel{\eqref{ip:chain}}{ \le} \frac{1}{\mu^\frac{1}{p-1}}\rho(u_n)\left(\frac{d(u_n)}{z(u_n)}\right)^\frac{1}{p-1}= \frac{1}{\mu^\frac{1}{p-1}}(c^*g(u_n)-f(u_n)-\dot z(u_n))\stackrel{\eqref{eq:seq*}}{<}1. \end{align*}
Hence,  there exists a value \(\varepsilon<\rho\) such that 
 \(z(\varepsilon)>\xi(\varepsilon)\). So, being \(\xi\) an upper-solution in $(0,\varepsilon]$, by Lemma \ref{lem:confronto} we infer \(z(u)\ge \xi(u) \) for every \(u\in (0,\varepsilon)\). Therefore, we 
   deduce that  
\[ \left(\frac{d(u)}{z(u)} \right)^\frac{1}{p-1}\le  \left(\frac{d(u)}{\xi(u)} \right)^\frac{1}{p-1} \stackrel{\eqref{ip:chain}}{<}M^\frac{1}{p-1}\frac{1}{\rho(u)} \quad \text{ for every } u\in (0,\varepsilon) \]
and this concludes the proof.
\end{pf}

\medskip
\begin{rem} \rm
In light of the previous theorem, the case in which $\ell_0=c^*g(0)-f(0)=0$  and $\dis\int_0^\frac12 \rho(u) \d u =+\infty$ remains unsolved. In this special case the solution can be smooth or not, as showed in \cite[Remark 10.1]{bcm2} in the particular case \(g(u)\equiv 1\).

 However, note that by virtue of Remark \ref{rem:1}, if $g(u)\equiv 1$ and $f(u)\equiv 0$ (i.e. a simple reaction-diffusion equation, without advection term), then necessarily $c^*>0$ and the previous theorem can be applied, without assuming \eqref{ip:chain2} (see the first two items of the statement).

\end{rem}

We conclude the section with a complete characterization of the slope of the t.w.s. when approaching the equilibrium 0. In what follows, \(\dot d(0)\) will denote the limit 
\(\dis\lim_{u\to 0} \frac{d(u)}{u}\), even if it is infinite.

\begin{thm}
\label{th:classif0}   Let   \eqref{ip:l0} be satisfied, with $\ell_0<+\infty$. 
Let \(u\) be a solution to \eqref{odeprobold} for some \(c\ge c^*\), defined in the interval \(I_u=(\alpha,\beta)\).

Then,  the following statements are true:

\smallskip
- \  if \(\dis \liminf_{u\to 0}d(u)> 0\), then \(	\dis \lim_{t \to \beta^-} u'(t)=0\), whatever \(c\ge c^*\) may be;

\smallskip
- \ if \(d(0^-)=0\), \ \(\dot d(0)=+\infty\){{,}} then \(\dis \lim_{t \to \beta^-} u'(t)=0\),   whatever \(c\ge c^*\) may be;

\smallskip
- \ if \(d(0^-)=0\), \ \(0<\dot d(0)<+\infty\), then \(\dis \lim_{t \to \beta^-} u'(t)=\begin{cases} 
-\left(\frac{c^*g(0)-f(0)}{\dot d(0)}\right)^\frac{1}{p-1} & \text{if } c=c^* \\ 0 & \text{if } c>c^*   \end{cases}\) 

\smallskip
- \ if \(d(0^-)=0\),  \(\dis\lim_{u\to 0}\dot d(u)= 0\)  and \ \(cg(0)-f(0)\ne 0\), then 
\[\dis \lim_{t \to \beta^-} u'(t)=\begin{cases}  -\infty & \text{ if } c=c^* \\ 0 & \text{ if } c>c^*.  \end{cases} \]

\end{thm}

\begin{pf} First of all, notice that assumption \eqref{ip:l0} guarantees  the validity of \eqref{ip:h0} for the function \(h(u):=(d(u))^\frac{1}{p-1} \rho(u)\).
Let \(u\) be a t.w.s for some \(c\ge c^*\). Since \(u'(t)=-\left(\dfrac{z(u(t))}{d(u(t))}\right)^\frac{1}{p-1}\) (see Proposition \ref{prop:equiv}), we have \(\dis\lim_{t\to b^-}u'(t)=-\lim_{u\to 0} \left(\frac{z(u)}{d(u)}\right)^\frac{1}{p-1}\). So, if \(\dis\liminf_{u\to 0}d(u)\ne 0\), we have \(\dis\lim_{t\to b^-}u'(t)=0\). Similarly, if 
 \(d(0^-)=0\) and \(\dot d(0)=+\infty\), since \(z\) is differentiable at 0 by Lemma \ref{l:general}, we have \(\dis\lim_{t\to b^-}u'(t)=0\).
Therefore, the first two statements are true.

\smallskip
Assume now \(d(0^-)=0\) and \(\dot d(0)< +\infty\). In this case, we have \(h_0=0\) (see \eqref{ip:h0}). So, since \(cg(0)-f(0)\ge 0\)  by \eqref{eq:stimanuova} and $\eta_0(0)=0$, by applying Lemma 
\ref{lem:zderiv0} 
we get 
\[ \dot z(0)= \begin{cases}  c^*g(0)-f(0) & \text{ if  } c=c^* \\ 0 & \text{ if } c>c^*. \end{cases}\]
and then the third item of the assertion  follows.

Assume now 
 \(d(0^-)=\dot d(0)= 0\), \(cg(0)-f(0)\ne 0\) and 
\(c=c^*\). Then, by Lemma \ref{lem:zderiv0} we have   \(\dot z(0)= c^*g(0)-f(0)>0\), hence  
\(\dis\lim_{t\to b^-}u'(t) = -\infty\) since \(\dot d(0)=0\). 

Consider now the last case, when  \(d(0^-)=0\), \(\dot d(u)\to 0\) as \(u\to 0\), \(cg(0)-f(0)\ne 0\) and  \(c>c^*\).
Since \(\dot z(0)=0\), there exists a sequence \(u_n\to 0\) such that \(\dot z(u_n)\to 0\). Therefore, by  \eqref{odeprob2}  we get 
\beq
\left(\frac{d(u_n)}{z(u_n)}\right)^\frac{1}{p-1} = \frac{cg(u_n)-f(u_n)-\dot z(u_n)}{\rho(u_n)}\to +\infty \label{e:succ}
\eeq
since \(cg(0)-f(0)>0\).

Let us now fix a real \(\varepsilon>0\) and define 
\(\varphi_\varepsilon(u):=\varepsilon d(u)\).
Notice that 
\[ \lim_{u\to 0^+} cg(u)-f(u)-\frac{\rho(u)(d(u))^\frac{1}{p-1}}{(\varphi_\varepsilon(u))^\frac{1}{p-1}}= cg(0)-f(0)>0.\]
So, since \(\dot \varphi_\varepsilon(u)=\varepsilon \dot d(u)\to 0\) as \(u\to 0\), there exists a value \(\delta=\delta_\varepsilon>0\) such that
\[
\dot \varphi_\varepsilon(u) < cg(u)-f(u)-\frac{\rho(u)(d(u))^\frac{1}{p-1}}{(\varphi_\varepsilon(u))^\frac{1}{p-1}} \quad \text{ for every } u\in (0,\delta),
\]
that is \(\varphi_\varepsilon\) is a  lower-solution for the equation of \eqref{odeprob2} in the interval \((0,\delta)\).

By \eqref{e:succ} we have \( \frac{d(u_n)}{z(u_n)}>\frac{1}{\varepsilon}\) for \(n\) sufficiently large; so there exists a value \(\eta=\eta_\varepsilon< \delta\) such that \(z(\eta)< \varphi_\varepsilon(\eta)\).
Since \(\varphi_\varepsilon\) is a  lower-solution, by Lemma \ref{lem:confronto} we derive \(z(u)< \varphi_\varepsilon(u)=\varepsilon d(u)\) for every \(u\in (0, \eta)\). Then, 
\( 0< \frac{z(u)}{d(u)}< \varepsilon \quad \text{ for every } u\in (0,\eta)\),
 implying that \(\dis\lim_{t\to b^-} u'(t)=0\) and this concludes the proof.
\end{pf}

\bigskip

\begin{rem}\rm
Both in Theorem \ref{t:T2*} and \ref{th:classif0} the value of $c^*g(0)-f(0)$ plays a relevant role. Recall that Theorem \ref{t:stima} (see also Remark \ref{rem:c*}) provides a sufficient condition and a necessary one in order to establish if $c^*g(0)-f(0)$ is positive or null. 

Moreover, note that the results of Theorem \ref{th:classif0} agree with those of Theorem  \ref{t:T2*}. In fact, when  $d(0)=0$ and $\dot d(0)<+\infty$, but $c^*g(0)\ne f(0)$ (that is the last two items of Theorem \ref{th:classif0}), then the integral 
$\dis\int_0^\frac12 \left(\frac{d(u)}{u}\right)^\frac{1}{p-1}\!\! \d u$ converges, so $\beta<+\infty$, in accordance with   $\dis\lim_{t\to b^-} u'(t)<0$.
\end{rem}

\medskip

As a consequence of the results above proved, we can state the following simple criterium applicable when $d$ and $\phi$ are asymptotic to a power of $u$.

\begin{cor}\label{cor:0}
Assume that $d(u)\sim k_1 u^\delta$ and $\rho(u)\sim k_2 u^r$ as $u \to 0$, for some positive constants $k_1, k_2,r$ and $\delta \in \R$, satisfying 
\beq \label{ip:nec0}
r\ge \frac{1- \delta}{p-1}. \eeq 
Let $u$ be the t.w.s. having speed $c^*$.
Then,
\begin{itemize}
\item \ if \ $r<1$, we have $\beta<+\infty$;
\item \ if \ $r=\frac{1-\delta}{p-1} \ge 1$, we have $\beta=+\infty$;
\item \ if \ $r> \frac{1-\delta}{p-1}$ and $c^*g(0)>f(0)$, we have $\beta<+\infty$ if and only if \ $\frac{1-\delta}{p-1}<1$.

\end{itemize}

Moreover, when $\beta<+\infty$ we have
\begin{itemize}
\item[$\ast$] \ if $\delta<1 $, then $u'(\beta^-)=0$;
\item[$\ast$] \ if $\delta=1$  then $u'(\beta^-)=0$ if and only if 
$c^*g(0)-f(0)=0$;
\item[$\ast$] \  if $\delta> 1$ and $c^*g(0)>f(0)$, then $u'(\beta^-)=-\infty$.
\end{itemize}

\begin{pf} The proof consists is an immediate application of the results stated in this section. We just observe that condition \eqref{ip:nec0} is equivalent to the assumption $\ell_0<+\infty$, necessary in order to have t.w.s.
\end{pf}

\end{cor}

\begin{rem} \label{r:confr}\rm
The previous result extends \cite[Theorem 5.3]{DJKZ}, since for $c=c^*$ this only consider that case $0<r<1$.
\end{rem}

\section{Local analysis at the equilibrium 1}\label{sec:1}

In this section we study the property of {\em finite speed of saturation} of equation \eqref{eq:RD}, occurring when the solution reaches the equilibrium 1 in a finite time, that is the extremum $\alpha$ of the existence interval  is finite.

\medskip
First  of all, observe that 
\beq \label{eq:integr2}
\alpha= \int_{\alpha}^0  \left(\frac{d(u(t))}{z(u(t))} \right)^\frac{1}{p-1} \!\! u'(t)\d t=- \int_\frac12^1 \left(\frac{d(u)}{z(u)}\right)^\frac{1}{p-1}\!\!  \d u\eeq
So, \(\alpha\) is finite or infinite according to the convergence or the divergence of the last integral in \eqref{eq:integr2}.

\medskip
In what follows, we will assume that
\beq \label{ip:l1}  \text{ there exists the limit } \ \ell_{1}:= \lim_{u\to 1^-}\rho(u)\left(\frac{d(u)}{(1-u)}\right)^\frac{1}{p-1}\in [0,+\infty].\eeq
The first criterium concerns the case when $\ell_1<+\infty$.

\medskip

\begin{thm}
\label{t:T1}  Assume    \eqref{ip:l1} is satisfied, with $\ell_1<+\infty$.
Let \(u\) be a solution of \eqref{odeprobold} for some \(c\ge  c^*\).
Then the following statements hold true:

\smallskip {\bf \rm (a)} \
 if $\ell_1>0$    we have 
 
\[ \alpha>-\infty\quad \text{ if and only if } \quad \dis\int_\frac12^1 \frac{1}{\rho(u)} \d u<+\infty;\]

\smallskip {\bf \rm (b)} \
 if $\ell_1=0$  and \ $cg(1)-f(1)<0$,  we have 
 
\[ \alpha>-\infty\quad \text{ if and only if } \quad \dis\int_\frac12^1 \!\!\left( \frac{d(u)}{1-u}\right)^\frac{1}{p-1}\!\! \d u<+\infty;\]

\smallskip {\bf \rm (c)} \
 if $\ell_1=0$, \ $ cg(1)-f(1)>0$  and 
 there exist three constants \(\delta, \mu, M>0\) and a  function 
 $\xi\in C^1([1-\delta,1])$  such that  
\beq \label{ip:chain2}
\mu \xi(u) \le  d(u)(\rho(u))^{p-1}\le M \xi(u) \quad \text{ for every } u\in [1-\delta,1), 
\eeq
we have
\[\alpha>-\infty \quad \text{ if and only if } \quad \dis\int_\frac12^1 \frac{1}{\rho(u)} \d u<+\infty;\]

\smallskip {\bf \rm (d)} \
 if $\ell_1= cg(1)-f(1) =0$ and \eqref{ip:chain2} holds true, then
 \[ \dis \int_\frac12^1 \frac{1}{\rho(u)}\d u<+\infty \quad \text{ implies } \quad \alpha>-\infty.\]  

\end{thm}
 
\begin{pf}
 Let $z$ be the  solution of problem \eqref{odeprob2} corresponding to the t.w.s. $u$ (see Proposition \ref{prop:equiv}). Set $h(u):=\rho(u)(d(u))^\frac{1}{p-1}$.
 
\medskip 
\(\bullet\) {\bf Statement $\boldsymbol{(a)}$:} \ \(\ell_1>0\).

\smallskip
By virtue of Lemma \ref{lem:zderiv1}  we have  $\dot z(1)<0$. 
Therefore,
	\[ \frac{\ell_1}{\rho(u)}\sim  \left(\frac{d(u)}{(1-u)}\right)^\frac{1}{p-1}\sim |\dot z(1)|^\frac{1}{p-1} \left( \frac{d(u)}{z(u)}\right)^\frac{1}{p-1} \quad \text{as } u \to 1\]
and the assertion follows from \eqref{eq:integr2}.

\medskip 
\(\bullet\) {\bf Statement $\boldsymbol{(b)}$:} \ \(\ell_1=0\) and $cg(1)-f(1)<0$.

\smallskip 
Also in this case
by virtue of Lemma \ref{lem:zderiv1}  we have  $\dot z(1)<0$. 
Therefore, 
	\[  \frac{d(u)}{z(u)}\sim   \frac{d(u)}{|\dot z(1)|(1-u)} \quad \text{as } u \to 1\]
and the assertion follows from \eqref{eq:integr2}.

\medskip

\(\bullet\) {\bf Statement $\boldsymbol{(c)}$:} \ \(\ell_1=0\) and $cg(1)-f(1)>0$, under assumption \eqref{ip:chain2}.
\smallskip

 Put 
\[ \lambda_1:=cg(1)-f(1)> 0, \quad
k_1:=\left(\frac{2}{\lambda_1}\right)^{p-1} \ , \quad
\psi(u):=Mk_1\xi(u).\]
Of course, \(\psi\in C^1([1-\delta,1])\), with 
$\dot \psi(1)=0$ by \eqref{ip:chain2}, since $\ell_1=0$. Moreover, again by \eqref{ip:chain2} we have that 
\[ \limsup_{u \to 1} \rho(u) \left(\frac{d(u)}{\psi(u)}\right)^\frac{1}{p-1}\le \frac{1}{k_1^\frac{1}{p-1}}=\frac12 \lambda_1.\]  
Therefore, we have
\begin{align*}  \dot\psi(1) & =0 < \frac12\lambda_1 \le
\liminf_{u\to 1} \left(cg(u)-f(u) - \frac{(d(u))^\frac{1}{p-1} \rho(u)}{(\psi(u))^\frac{1}{p-1}}\right).\end{align*}
Then, for some \(\rho\in(0,\delta)\) we have 
\[ \dot\psi(u) < cg(u)-f(u) - \frac{(d(u))^\frac{1}{p-1} \rho(u)}{(\psi(u))^\frac{1}{p-1}} \quad \text{ for every } u \in (1-\rho,1)\]
i.e. \(\psi\) is a lower-solution in the interval \((1-\rho,1)\). Hence, by Lemma \ref{l:comp1} we deduce $z(u)\le \psi(u)$ for every $u\in (1-\delta,1)$.

Therefore,  we obtain
\beq \left(\frac{d(u)}{z(u)} \right)^\frac{1}{p-1}\ge  \left(\frac{d(u)}{\psi(u)} \right)^\frac{1}{p-1} \stackrel{\eqref{ip:chain2}}{>}\left(\frac{\mu}{Mk_1}\right)^\frac{1}{p-1}\!\!\frac{1}{\rho(u)} \quad \text{ for every } u\in (1-\rho,1). \label{eq:stimasotto2} \eeq

\medskip

Put now  $\varphi(u):=\dfrac{1}{3^{p-1}}\mu k_1\xi(u)$.
Notice that 
\[\dot \varphi(1)=0>\lambda_1-\frac{3}{k_1^\frac{1}{p-1}}\ge
\limsup_{u \to 1} \left(cg(u)-f(u)- \frac{(d(u))^\frac{1}{p-1}\rho(u)}{(\varphi(u))^\frac{1}{p-1}}\right).
\]
Therefore, for some $\sigma\in (0,\rho)$ we have
\[\dot \varphi(u)> cg(u)-f(u)- \frac{(d(u))^\frac{1}{p-1}\rho(u)}{(\varphi(u))^\frac{1}{p-1}} \quad \text{ for every } u\in (1-\sigma,1),\]
that is $\varphi$ is an upper-solution in $(1-\sigma,1)$.
Then, again by virtue of Lemma \ref{l:comp1} we infer
 that  $\varphi(u) \le z(u)$ in $(1-\sigma,1)$.
Thus, 
\beq\left(\frac{d(u)}{z(u)} \right)^\frac{1}{p-1}\le
\frac{3}{(\mu k_1)^\frac{1}{p-1}}
 \left(\frac{d(u)}{\xi(u)} \right)^\frac{1}{p-1} \stackrel{\eqref{ip:chain2}}{\le }
 3\left(\frac{M}{\mu k_1}\right)^\frac{1}{p-1}
 \frac{1}{\rho(u)} \quad \text{in } (1- \sigma,1). \label{eq:stimasopra2} \eeq
Therefore, the assertion follows from \eqref{eq:stimasotto2} and \eqref{eq:stimasopra2}.

\medskip {\bf Statement $\boldsymbol{(d)}$:} \ \(\ell_1=0\) and $cg(1)-f(1)=0$, under assumption \eqref{ip:chain2}.

\smallskip

\noindent
Note that 
$\dot \xi(1)=0$ by \eqref{ip:chain2}, since $\ell_1=0$. Thus, again by \eqref{ip:chain2} we have 
\begin{align*}  \dot\xi(1) & =0 > -\mu^\frac{1}{p-1} \ge
\limsup_{u\to 1} \left(cg(u)-f(u) - \frac{(d(u))^\frac{1}{p-1} \rho(u)}{(\xi(u))^\frac{1}{p-1}}\right).\end{align*}
Then, for some \(\rho\in(0,\delta)\) we have 
\[ \dot\xi(u) > cg(u)-f(u) - \frac{(d(u))^\frac{1}{p-1} \rho(u)}{(\xi(u))^\frac{1}{p-1}} \quad \text{ for every } u \in (1-\rho,1)\]
i.e. \(\xi\) is an  upper-solution in the interval \((1-\rho,1)\).
 So, we again infer from Lemma \ref{l:comp1} that \(\xi(u)\le z(u)\) for every \(u\in (1-\rho,1)\). Therefore, we 
   deduce that  
\[ \left(\frac{d(u)}{z(u)} \right)^\frac{1}{p-1}\le  \left(\frac{d(u)}{\xi(u)} \right)^\frac{1}{p-1} \stackrel{\eqref{ip:chain2}}{<}M^\frac{1}{p-1}\frac{1}{\rho(u)} \quad \text{ for every } u\in (1-\rho,1) \]
and this concludes the proof.
\end{pf} 

\medskip
\begin{rem} \rm
In light of Theorem \ref{t:T1}, the case $\ell_1=cg(1)-f(1)=0$ and $\dis\int_\frac12^1 \frac{1}{\rho(u)}\d u=+\infty$ remains unsolved. However, note that by Remark \ref{rem:1} when $g(u)\equiv 1$ and $f(u)\equiv 0$ (that is when we deal with a mere reaction-diffusion equation without convection effects) we have $c>0$, so Theorem \ref{t:T1} can be applied (see statement $(c)$). 

\end{rem}

\medskip
The next result covers the cases when $\ell_1=+\infty$.

\medskip

\begin{thm}\label{t:infinity}
Suppose that for some positive constants $c_1, c_2, \delta$ and $\lambda\in (-1, \frac{1}{p-1}]$ we have
\beq\label{ip:asint}
c_1 (1-u)^\lambda \le \rho(u)(d(u))^\frac{1}{p-1} \le c_2 (1-u)^\lambda \quad \text{ in } (1-\delta,1).
\eeq
Then, 
\beq \label{eq:integral} \alpha = -\infty \ \ \text{ if and only if } \ \ \int_0^\frac12 \frac{(1-u)^{\lambda- \frac{\lambda+1}{p}}}{\rho(u)} \d u =+\infty.\eeq
Moreover, when \ $\alpha>-\infty$, then $\dis\lim_{t \to \alpha^+}u'(t)=0$.

\end{thm}

\begin{pf}
Put $M:=\dis\max_{u\in [0,1]}(cg(u)-f(u))$, let $k>0$ be such that 
\beq  \label{eq:K1}\frac{c_1}{k^\frac{1}{p-1}}-k(\lambda+1)\frac{p-1}{p}>M\eeq
and put $\psi(u):=k (1-u)^{(\lambda+1)\frac{p-1}{p}}$, $u\in [0,1]$.
We have $\psi(1)=0$ and $\psi(u)>0$ in $(0,1)$. 
Moreover, since $\lambda\le \frac{1}{p-1}$, we have $\lambda\le \frac{\lambda+1}{p}$, so 
\beq \label{eq:eq2}\left(\frac{c_1}{k^\frac{1}{p-1}}-k(\lambda+1)\frac{p-1}{p} \right)(1-u)^{\lambda- \frac{\lambda+1}{p}} \ge \frac{c_1}{k^\frac{1}{p-1}}-k(\lambda+1)\frac{p-1}{p}\stackrel{\eqref{eq:K1}}{>} cg(u)-f(u)\eeq
for every $u\in (1-\delta,1)$. Then, for every $u\in (1-\delta,1)$ we have
\begin{align*} \dot \psi(u)= & -k(\lambda+1)\frac{p-1}{p} (1-u)^{\lambda- \frac{\lambda+1}{p}} \\ \stackrel{\eqref{eq:eq2}}{>} & \ cg(u)-f(u) - \frac{c_1}{k^\frac{1}{p-1}} (1-u)^{\lambda- \frac{\lambda+1}{p}} \\
\stackrel{\eqref{ip:asint}}{>}  & \ cg(u)-f(u)- \frac{ \rho(u)(d(u))^\frac{1}{p-1}}{(\psi(u))^\frac{1}{p-1}}. 
\end{align*}

\noindent
So,  $\psi$ is an upper-solution to the equation in \eqref{odeprob2} in $(1-\delta,1)$. Thus, by virtue of Lemma \ref{l:comp1} we obtain $\psi(u)\le z(u)$ in $(1-\delta,1)$, yielding 
\[ \left(\frac{d(u)}{z(u)}\right)^\frac{1}{p-1}  \le \left( \frac{d(u)}{\psi(u)}\right)^\frac{1}{p-1} \stackrel{\eqref{ip:asint}}{\le}   \frac{c_2(1-u)^\lambda}{\rho(u)(\psi(u))^\frac{1}{p-1}} = \frac{c_2}{k^\frac{1}{p-1}} \frac{(1-u)^{\lambda-\frac{\lambda+1}{p}}}{\rho(u)}.\]
Therefore, if the integral in \eqref{eq:integral} converges, by \eqref{eq:integr2}  we deduce that $\alpha>-\infty$.

\bigskip
In order to prove the converse, put $m:=\dis\min_{u\in [0,1]}(cg(u)-f(u))$, let $\ell>0$ be such that 
\[\frac{c_2}{\ell^\frac{1}{p-1}}-\ell(\lambda+1)\frac{p-1}{p}<\min\{m,0\}\]
and put $\varphi(u):=\ell (1-u)^{(\lambda+1)\frac{p-1}{p}}$, $u\in [0,1]$.
We have $\varphi(1)=0$ and $\varphi(u)>0$ in $(0,1)$. 
Moreover, 
similarly to what we have done above, it is possible to show that $\varphi$ is a lower-solution to \eqref{odeprob2} in $(1-\delta,1)$. Thus, by virtue of Lemma \ref{l:comp1} we obtain $\varphi(u)\ge z(u)$ in $(1-\delta,1)$, yielding
 
\beq \label{eq:dis1} \left(\frac{d(u)}{z(u)}\right)^\frac{1}{p-1}  \ge \left( \frac{d(u)}{\varphi(u)}\right)^\frac{1}{p-1} \stackrel{\eqref{ip:asint}}{\ge}   \frac{c_1(1-u)^\lambda}{\rho(u)(\varphi(u))^\frac{1}{p-1}} = \frac{c_1}{\ell^\frac{1}{p-1}} \frac{(1-u)^{\lambda-\frac{\lambda+1}{p}}}{\rho(u)}.\eeq
Therefore, if the integral in \eqref{eq:integral} diverges, we deduce that $\alpha=-\infty$ and this concludes the proof of \eqref{eq:integral}. 

Finally, by  \eqref{eq:dis1} we have 
\[  \left(\frac{z(u)}{d(u)} \right)^\frac{1}{p-1}\le \frac{\ell^\frac{1}{p-1}}{c_1} \frac{\rho(u)}{(1-u)^{\lambda-\frac{\lambda+1}{p}}} \quad \text{ in } \ (1-\delta,1)\]
Therefore, since $\lambda\le \frac{\lambda+1}{p}$, we have $\dis\lim_{u\to 1^-} \left(\frac{z(u)}{d(u)} \right)^\frac{1}{p-1}=0$ and the assertion follows from the equation in  \eqref{pr:cau}.
\end{pf}

\bigskip

We conclude this investigation with a simple sufficient condition  for the unboundedness of $\alpha$ when $\ell_1=+\infty$, which can be applied when the function $\rho d^\frac{1}{p-1}$ is not asymptotic to a power of $(1-u)$ as $u\to 1$.

\medskip
\begin{thm}\label{th:little}
 Put $\varphi(u):=d(u) (\rho(u))^{p-1}$, suppose that $\varphi$ is differentiable in a left neighborhood of 1 and satisfies \[ \varphi(1^-)=0 \quad \text{ and } \quad \dot \varphi(u) \to -\infty \text { as } u\to 1^-.\]
Then, if $\dis\int_\frac12^1 \frac{1}{\rho(u)} \d u =+\infty$, we have $\alpha=-\infty$.
\end{thm}
\begin{pf}
By the assumptions, we have that $\varphi$ is a positive lower-solution to the equation in \eqref{odeprob2}  in a left neighborhood of 1, then by  Lemma \ref{l:comp1} we obtain $\varphi(u)\ge z(u)$ in $(1-\delta,1)$, yielding 
\[ \left(\frac{d(u)}{z(u)}\right)^\frac{1}{p-1}  \ge \left( \frac{d(u)}{\varphi(u)}\right)^\frac{1}{p-1} = \frac{1}{\rho(u)}\]
and the assertion follows from \eqref{eq:integr2}.
\end{pf}

\bigskip

We state now the complete study of the slope of the t.w.s. when approaching the   equilibrium 1.  In what follows, \(\dot d(1)\) will denote the limit \(\dis\lim_{u\to 1}\frac{-d(u)}{1-u}\), even if it is infinite.

\begin{thm}
\label{th:classif1} Assume that  \eqref{ip:l1} holds true.
Let \(u\) be a solution of \eqref{odeprobold}, defined in the interval \(I_u=(\alpha,\beta)\).

Then, the following statements are true:

\smallskip - \
 if \(\dis\liminf_{u \to 1}d(1)> 0\), then \(	\dis \lim_{t \to \alpha^+} u'(t)=0\);
 
\smallskip - \ if \(d(1^-)=0\), \ \( \dot d(1)=-\infty\) \ and \ $\ell_1<+\infty$, then \(\dis \lim_{t \to  \alpha^+} u'(t)=0 \);

\smallskip - \ if \(d(1^-)=0\), \  \(-\infty < \dot d(1)< 0\){{,}} then \(\dis \lim_{t \to  \alpha^+} u'(t)=-\left(
\max\left\{0, \frac{cg(1)-f(1)}{\dot d(1)}\right\}\right)^{\frac{1}{p-1}} \);

\smallskip - \ if \(d(1^-)=0\),  \(\dis\lim_{u\to 1}\dot d(u)= 0\)  and \(cg(1)-f(1)\ne 0\), then 
\[\dis \lim_{t \to  \alpha^+} u'(t)=\begin{cases}  0 & \text{ if } \ cg(1)-f(1)>0  \\ -\infty & \text{ if }\ cg(1)-f(1)<0.  \end{cases} \]

\end{thm}

\begin{pf}
First of all, notice that assumption \eqref{ip:l1} guarantees  the validity of \eqref{ip:h1} for the function \(h(u):=(d(u))^\frac{1}{p-1} \rho(u)\).

Let \(u\) be a t.w.s for some \(c\ge c^*\).
Since \(u'(t)=-\left(\dfrac{z(u(t))}{d(u(t))}\right)^\frac{1}{p-1}\) (see Proposition \ref{prop:equiv}), we have \(\dis\lim_{t\to \alpha^+}u'(t)=-\lim_{u\to 1^-} \left(\frac{z(u)}{d(u)}\right)^\frac{1}{p-1}\). So, if \(\dis\liminf_{u \to 1}d(1)> 0\), then \(\dis\lim_{t\to \alpha^+}u'(t)=0\).  Similarly, if
 \(d(1^-)=0\), \(\dot d(1)=-\infty\) and $\ell_1<+\infty$, since \(z\) is differentiable at 1 by Lemma \ref{l:general}, we have \(\dis\lim_{t\to a^+}u'(t)=0\). 
Therefore, the first two statements are true.

\smallskip
Assume now \(d(1^-)=0\) and \(-\infty < \dot d(1)< 0\). Then,  \(h_1=0\) and  by Lemma \ref{lem:zderiv1} we have
\( \dot z(1)= \min\{0, cg(1)-f(1)\}. \)
Hence, for every \(c\ge c^*\) we have
\[\dis \lim_{t\to \alpha^+} u'(t)=\lim_{u\to 1^-}-\left( \frac{z(u)}{d(u)}\right)^\frac{1}{p-1}= -\left(\frac{\dot z(1)}{\dot d(1)}\right)^\frac{1}{p-1}=-\left(\max\left\{0, \frac{cg(1)-f(1)}{\dot d(1)}\right\}\right)^\frac{1}{p-1}.\]

Assume now  \(d(1^-)=\dot d(1)= 0\) and \(cg(1)-f(1)< 0\). 
Since \(\dot z(1)= cg(1)-f(1)<0\) and \(\dot d(1)=0\), we have  \(\dis\lim_{u\to 1^-} \frac{z(u)}{d(u)}= +\infty\). So, \(\dis\lim_{t\to a^+} u'(t)=-\infty\).

Finally, assume  \(d(1^-)= 0\), \(\dot d(u)\to 0\) as \(u \to 1\) and \(cg(1)-f(1)> 0\). Consider the function 
\(\varphi_\varepsilon=\varepsilon d(u)\). 
Since \(\dot \varphi_\varepsilon(u)= \varepsilon\dot d(u)\to 0\) as \(u\to 1\),  for some \(\delta=\delta_\varepsilon>0\) we have
\beq \dot \varphi_\varepsilon (u) < cg(u)-f(u) - \frac{(d(u))^\frac{1}{p-1}\rho(u)}{(\varphi_\varepsilon(u))^\frac{1}{p-1}} \quad \text{for every } u\in (1-\delta,1), \label{eq:domin1}\eeq
that is \(\varphi_\varepsilon\) is a  lower-solution in \((1-\delta,1)\), with $\varphi_\varepsilon(1^-)=0$. Hence, by virtue of Lemma \ref{l:comp1} we have
 \(z(u)\le \varphi_\varepsilon(u)=\varepsilon d(u)\) for every  \(u\in (1-\delta,1)\), implying
 \(\dis \lim_{u \to 1} \frac{z(u)}{d(u)}=0\); so   \(\dis\lim_{t\to \alpha^+} u'(t)=0\).

\end{pf}

\begin{rem}
\rm
Note that the results of Theorem  \ref{th:classif1} agree with those of Theorem \ref{t:T2}. In fact, when  $d(1^-)=0$ and $\dot d(1)>-\infty$, but $cg(1)- f(1)<0$, then we have $\dis\lim_{t \to \alpha^+} u'(t) <0$, so $\alpha>-\infty$, according to the first statement of Theorem \ref{t:T1}, since the ratio $d(u)/(1-u)$  is bounded.

Moreover, we recall that Theorem \ref{t:stima} provides a sufficient condition and a necessary one in order to establish if $c^*> k$ or $c^*\le k$ for each costant $k\in \R$. Such criteria can be used by putting $k=f(1)/g(1)$ in order to establish if $c^*>f(1)/g(1)$ or $c^*\le f(1)/g(1)$.  

\end{rem}

Now we state a simple criterium useful in the case when $d$ and $\rho$ go to 0 as $u\to 1$ like a power of $(1-u)$. It is an immediate consequence of the previous results.

\begin{cor}\label{cor:1}
Assume that $d(u)\sim k_1 (1-u)^\delta$ and $\rho(u)\sim k_2 (1-u)^r$ as $u \to 1$, for some positive constants $k_1, k_2,r$ and $\delta \in \R$, satisfying 
\beq \label{eq:nec1}r(p-1)+\delta+p> 1. \eeq 
Let $u$ be a t.w.s. having speed $c$.

Then,
\begin{itemize}
\item[(a)] \ if \ $r(p-1)+\delta = 1$, we have 
$\alpha>-\infty$ if and only if \ $p+\delta>2$;
\item[(b)] \ if \ $r(p-1)+\delta > 1$ and $cg(1)<f(1)$, we have 
$\alpha>-\infty$ if and only if \ $p+\delta>2$;
\item[(c)] \ if \ $r(p-1)+\delta > 1$ and $cg(1)>f(1)$, we have 
$\alpha>-\infty$ if and only if \  $r<1$;
\item[(d)] \ if \ $r(p-1)+\delta > 1$, $cg(1)=f(1)$ and $r<1$, we have 
$\alpha>-\infty$.
\item[(e)] \ if \ $r(p-1)+\delta < 1$, we have $\alpha>-\infty$ if and only if \ $p+\delta>r+1$. 
\end{itemize}

Moreover, when $\alpha>-\infty$ we have
\begin{itemize}
\item[$\ast$] \ if $\delta<1 $, then \ $u'(\alpha^-)=0$;
\item[$\ast$] \ if $\delta=1$,  then $u'(\alpha^-)=0$ if and only if 
\ $cg(1)-f(1)\ge 0$;
\item[$\ast$] \  if $\delta> 1$ and $cg(1)>f(1)$, then \ $u'(\alpha^-)=0$
\item[$\ast$] \  if $\delta> 1$ and $cg(1)<f(1)$, then \ $u'(\alpha^-)=-\infty$.
\end{itemize}

\begin{pf} Note that assumption \eqref{eq:nec1} ensures the validity of $(H_d)$. 
As for the statements (a)-(d), they are consequences of Theorem \ref{t:T1}. Instead, statement (e) follows from  Theorem \ref{t:infinity}.

Moreover, when $\alpha>-\infty$ and $\delta\le 0$ we have $u'(\alpha^-)=0$ as a consequence of Theorem \ref{th:classif1}; the same occurs when $0<\delta<1$ and $r(p-1)+\delta \ge 1$. Instead, when 
$0<\delta<1$ and $r(p-1)+\delta < 1$ the results follows from Theorem \ref{t:infinity}. Finally, the last assertions follow from Theorem \ref{th:classif1}.
\end{pf}

\end{cor}

\begin{rem}\rm  \label{rem:DR}
The study of the boundedness of $\alpha$ and the regularity of the solution $u$ was carried out also in \cite{DJKZ} (see also \cite{DZ} for the case without convection), in the case  $g(u)\equiv 1$ and the functions $d,\rho$ are asymptotic to powers of $(1-u)$ as $u\to 1$. It is possible to check that even in such a more particular setting the present results extend those obtained in \cite[Theorems 5.1, 5.2]{DJKZ}. In fact, in \cite{DJKZ} the authors assume $r(p-1)+\delta\ge 0$, while here we assume the more general condition \eqref{eq:nec1}.
Moreover, 
 we study also the case in which 
 $\ell_1=0$ and $c<f(1)$ (that is statement $(b)$ of Corollary \ref{cor:1}), which was not considered in \cite{DJKZ}.

Finally, let's mention that in \cite[Theorem 5.2]{DJKZ}
it is considered also the case when $\ell_1=0$ and $c-f(t)>0$ a.e. in a left neighborhood of 0, but as it is possible to check in light of the proof of that result, this works only if $c>f(1)$.

\end{rem}

\medskip

\section{Appendix}

We here discuss the validity of the results stated in Section \ref{sec:ODE}, taken from the paper \cite{Ma},
by explaining which parts of the proofs in \cite{Ma} must be modified in light of the present weaker assumption \eqref{ip:weakH} (instead of the previous condition $h\in C[0,1]$ with $h(0)=h(1)=0$ considered in \cite{Ma}).

\bigskip
{\bf Proof of Lemma \ref{l:general}.} 

\smallskip Note that $\dot z(u)<cg(u)-f(u)$ in $(u_0,u_1)$, so   $\dot z$ is bounded from above and then both the limits $z(u_0^+)$ and $z(u_1^-)$ exist and the latter is finite. Moreover, if $z(u_0^+)=+\infty$, then from \eqref{eq:sing} we derive $\dot z(u)\to cg(u_0)-f(u_0)$ as $u\to u_0^+$, a contradiction. Hence, $z(u_0^+)$ is finite too. 

Furthermore, if $u_0>0$ then necessarily $z(u_0^+)=0$, since $(u_0,u_1)$ is the maximal existence  interval of $z$. Hence, we have
$(cg(u)-h(u))(z(u))^\frac{1}{p-1}-h(u)\to -h(u_0)<0$  as $u\to u_0^+$, implying that 
\[ \dot z(u) =cg(u)-f(u)-\frac{h(u)}{(z(u))^\frac{1}{p-1}} < -\frac12 \frac{h(u_0)}{(z(u))^\frac{1}{p-1}} <0, \]
in a right neighborhood of $u_0$, a contradiction. So, $u_0=0$.

\medskip
Assume now $z(0^+)=0$ and \eqref{ip:h0}. Under these conditions we have $h_0<+\infty$. Indeed,  put $M:=\dis\max_{u\in [0,1]}|cg(u)-f(u)|$ we have  $0<z(u)\le \dis Mu$ in $[u_0,u_1]$, then
\[ \dot z(u)\le M - \frac{1}{M^\frac{1}{p-1}} \frac{h(u)}{u^\frac{1}{p-1}} \quad \text{ in } (0,u_1)\]
so, if $h_0=+\infty$ then  $\dot z(u)\to -\infty$ as $u\to 0^+$, a contradiction. Hence, $h_0<+\infty$ and then $h(0)=0$.

\smallskip
Moreover, it is possible to prove that there exists
 the limit $\lambda:=\dis\lim_{u\to 0} \frac{z(u)}{u}\in [0,+\infty]$ (the proof is exactly the same proposed in \cite[Lemma 1]{Ma}). 
Furthermore, if  $\lambda=+\infty$, then 
from
\beq \dot z(u)= cg(u)-f(u) - \frac{h(u)}{u^\frac{1}{p-1}}\frac{u^\frac{1}{p-1}}{(z(u))^\frac{1}{p-1}} \label{eq:eq}\eeq
we deduce that $\dot z(u)\to cg(0)-f(0)\le M$ as $u\to 0^+$, a contradiction. 
So $\lambda<+\infty$, implying that $z$ is differentiable at 0, with $\dot z(0)=\lambda$.

Now, if $\dot z(0)>0$, then from \eqref{eq:eq} we derive the existence of the limit $\dot z(0^+)$, which necessarily coincides with $\dot z(0)$ and passing to the limit as $u\to 0$ in \eqref{eq:eq} we infer  that $\eta_0(\dot z(0))=0$.
Instead, if $\dot z(0)=0$, then $h_0=0$ too, otherwise $\dot z(u)\to -\infty$, a contradiction. Therefore, in this case $\eta_0(\dot z(0))=\eta_0(0)=0$.

\medskip
Finally, assume $u_1=1$, $z(1^+)=0$ and \eqref{ip:h1}. In this case, it may happen that 
  $h_1=+\infty$. Nevertheless, if we assume $h_1<+\infty$, the proof of the existence of $\dot z(1)$, with $\eta_1(\dot z(1))=0$, is analogous to the previous one, related to the equilibrium 0. 
\hfill $\Box$

\bigskip
\begin{rem} \rm Note that  the statement of \cite[Lemma 1]{Ma} concerning the equilibrium 1 was not correct. For the equilibrium $1$ is not guaranteed that $h_1<+\infty$; in order to ensure the existence of $\dot z(1)$, one need to assume $h_1<+\infty$. 
\end{rem}

\bigskip
{\bf Proof of Proposition \ref{l:soprasol}.} 

\medskip
The proof given in \cite[Proposition 1]{Ma}  also works under the present weaker hypothesis \eqref{ip:weakH}, with the expedient described below.

In the case $\varphi(1^-)>0$,  even if $h$ is not continuous at 1 we can 
still consider 
the Cauchy problem given by the equation \eqref{eq:sing} and the condition  $z(1)=\varphi(1^-)/n>0$. Indeed, the classical Carath\'eodory conditions for the existence and uniqueness are satisfied, so there exists the solution $z_n$, which is absolutely continuous in $(0,1]$ and $C^1$ in $(0,1)$. Moreover, since $z_n(u)<\varphi (u)$ in a left neighborhood of $1$, we can again apply the comparison result stated in \cite[Lemma 2]{Ma} to deduce that $z_n(u)\le \varphi(u)$ in $(0,1]$ and the proof proceeds as in \cite{Ma}.

\hfill $\Box$

\bigskip
{\bf Proof of Theorem \ref{t:main}.} 

\medskip
In the proof of the analogous result given in \cite{Ma} (see Theorem 1), the continuity of $h(u)$ at $u=1$  is never used. The proof works even under the sole assumption that $h\in L^1(0,1)$. 
\hfill $\Box$.

\subsection*{Funding Information}
This article was supported by PRIN 2022 -- Progetti di Ricerca di rilevante Interesse Nazionale, ``Nonlinear differential problems with applications to real phenomena'' (2022ZXZTN2). 

\subsection*{Acknowledgements}
The author is member of the Gruppo Nazionale per l'Analisi Matematica, la Probabilit\`a e le loro Applicazioni (GNAMPA) of the Istituto Nazionale di Alta Matematica (INdAM).

\end{document}